\documentclass[leqno]{amsart}

\usepackage{amsthm,amssymb,amsmath,amscd,latexsym}

\usepackage[latin1]{inputenc}
\usepackage[all]{xy}

\theoremstyle{plain}
\newtheorem{theorem}{Theorem}[section]
\newtheorem{proposition}[theorem]{Proposition}
\newtheorem{lemma}[theorem]{Lemma}
\newtheorem{corollary}[theorem]{Corollary} 
 
\theoremstyle{definition}
\newtheorem{definition}[theorem]{Definition}
\newtheorem{remark}[theorem]{Remark}
\newtheorem{remarks}[theorem]{Remarks}
\newtheorem{example}[theorem]{Example}

\theoremstyle{remark}

\numberwithin{equation}{section}

\def\paritem#1{%
  \smallskip
  \setbox0=\hbox{#1\enspace}
  \par\noindent
  \ifnum\wd0>\parindent\box0
  \else\hbox to\parindent{\box0\hfill}\fi\ignorespaces}

\def\C{\mathbb{C}}
\def\F{\mathbb{F}}
\def\Z{\mathbb Z}
\def\R{\mathbb R}
\def\S{\mathbb S}
\def\PBbb{\mathbb P}
\def\g{\mathfrak g}
\def\D{\mathcal D}
\def\F{\mathcal F}
\def\I{\mathcal I}
\def\L{\mathcal L}
\def\P{\mathcal P}
\def\U{\mathcal U}
\def\V{\mathcal V}
\def\inv{^{-1}}
\def\Sing{{\rm Sing}}
\def\epsilon{\varepsilon}
\def\dd#1#2{\frac{d{#1}}{d{#2}}}

\DeclareMathOperator\Diff{Diff}
\DeclareMathOperator\Hom{Hom}
\DeclareMathOperator\id{id}
\DeclareMathOperator\SU{SU}

\newcommand{\what}{\widehat}
\newcommand{\wbar}{\overline}

\newcommand{\wti}{\widetilde}

\newcommand{\GL}{\operatorname{GL}}
\newcommand{\gl}{\operatorname{\mathfrak{gl}}}

\newcommand{\SO}{\operatorname{SO}}
\newcommand{\Un}{\operatorname{U}}

\newcommand{\Tr}{\operatorname{Tr}}

\newcommand{\Id}{\operatorname{Id}}
\newcommand{\pt}{\operatorname{pt}}

\newcommand{\ra}{\rightarrow}

\newcommand{\lra}{\longrightarrow}

\newcommand{\oset}[1]{\overset {#1}{\ra}}

\newcommand{\ocong}[1]{\overset {#1}{\cong}}

\newcommand{\x}{\times}

\newcommand{\ot}{\otimes}
\newcommand{\tp}[2]{{#1} \ot {#2}}

\newcommand{\Sky}{\mathbf{Sky}}

\newcommand{\Loc}{\mathbf{Loc}}

\newcommand{\veps}{\varepsilon}

\begin{document}
\title[Gerbes, simplicial forms and invariants \ldots]
{Gerbes, simplicial forms and invariants for families of foliated bundles}
\author[J.~L~Dupont]{Johan L.~Dupont${}^1$}
\address{Department of Mathematics   \\
University of Aarhus \\
DK-8000 {\AA}rhus C, Denmark}
\email[J.~L~Dupont]{dupont@imf.au.dk}
\author[F.~W~Kamber]{Franz ~W.~Kamber${}^2$}
\address{Department of Mathematics \\ 
University of Illinois at Urbana--Champaign \\
1409 W. Green Street \\
Urbana, IL 61801, USA}
\email[F.~W~Kamber]{kamber@math.uiuc.edu}

\thanks{
Work supported in part by the Erwin Schr\"odinger International 
Institute of Mathematical Physics, Wien, Austria and by the 
Statens Naturvidenskabelige Forskningsr{\aa}d, Denmark} 
\thanks{
${}^1$ Supported in part by the European Union Network EDGE. 
${}^2$ Supported in part by `Fonds zur F\"orderung der 
wissenschaftlichen  Forschung, Projekt P~14195~MAT'}
\date{\today}

\begin{abstract}
  The notion of smooth Deligne cohomology is conveniently reformulated
  in terms of the simplicial deRham complex. In particular the usual
  Chern-Weil and Chern-Simons theory is well adapted to this framework
  and rather easily gives rise to characteristic Deligne cohomology
  classes associated to families of bundles and connections. 
  In turn this gives invariants for families of foliated bundles.
  The construction provides representing cocycles in the usual
  \v Cech-deRham model for smooth Deligne cohomology called 
  `gerbes with connection' as they generalize usual Hermitian line 
  bundles with connection. A special case is the Quillen line bundle 
  associated to families of flat $\SU(2)$-bundles.
\end{abstract}

\maketitle

\tableofcontents



\medbreak
\section{Introduction}\label{one}

The determinant line bundle was constructed by Quillen~\cite{Q} for 
families of Riemann surfaces and generalized to higher dimension by
Bismut and Freed (see e.g. \cite{BF}, \cite{F1}, \cite{F2}). 
It also admits a `geometric' construction (and further generalization) 
in terms of families of principal $G$-bundles with connection for $G$ 
any Lie group (see e.g. Bonora et.al.~\cite{BCRS}, Brylinski~\cite{B2}, ~
\cite{B1}, Dupont--Johansen~\cite{DJ}). 

In this situation the construction in the present paper more generally 
provides `$\ell$--gerbes with connection' for suitable $\ell=0,1,2,\ldots$
depending on curvature conditions on the fibre connections in the family.
We use the phrase `(Hermitian line) gerbe' (respectively `(Hermitian
line) gerbe with connection') as an abreviation for the notion of a
representing cocycle (with a shift in degree) in the \v Cech 
(respectively \v Cech--deRham) model for the usual (respectively
Deligne) cohomology associated to the sheaf $\underline{U(1)}$ of
smooth functions with values in the circle group $U(1)\subseteq\C$.
We are aware that the word `gerbe' originally was used for a rather
different kind of object which however, in the abelian case, is closely
related to our `2--gerbe' in the same way as a `1--gerbe'
corresponds to a Hermitian line bundle. 
Similarly our notion of a `2--gerbe with connection' is in accordance 
with Hitchin~\cite{H} and is in line with a widespread use of the word 
`gerbe' in mathematical physics (see e.g. Carey--Mickelsson~\cite{CM}.) 
We refer to Brylinski~\cite{B2},~\cite{B3} for more information about 
Deligne cohomology and its relation to the original notion of `gerbes' 
(see also Breen--Messing~\cite{BM}). 
However, we are using the word `gerbe' only in the restricted sense
described in section 2.

Let us now describe our main results. 
In the following, $X$ will be a compact oriented smooth manifold 
and $G$ a Lie group with finitely many components.

\begin{definition}
  \label{1.1}
  A \emph{family} of principal $G$-bundles over $X$ with connections
  consists of the following:
  
  \paritem{(i)} A smooth fibre bundle $\pi\colon Y\to Z$ with fibre
  $X$ and structure group $\Diff^+(X)$ of orientation preserving
  diffeomorphisms.
  
  \paritem{(ii)} A principal $G$-bundle $p\colon E\to Y$.
  
  \paritem{(iii)} A smooth family $A=\{A_z\mid z\in Z\}$ of
  connections in the $G$-bundles $P_z = E~| {X_z}~, ~X_z = \pi\inv(z)$.
\end{definition}

Notice that the family of connections in (iii) can always be obtained
(using a partition of unity) from some `global' connection $B$ in
the $G$-bundle $E$ such that $A_z = B ~|~ T {P_z}$ for all $z\in Z$. 
But this \emph{global extension} is not part of the structure. 
Furthermore let $I_\Z^{n+1}(G)\subseteq I^{n+1}(G)$ denote the set of 
invariant homogeneous polynomials of degree $n+1$ on the Lie algebra 
$\g$ such that the Chern-Weil image is an \emph{integral} class. 
That is, $Q\in I_\Z^{n+1}(G)$ corresponds in the cohomology 
$H^{2n+2}(BG,\R)$, $BG$ the classifying space, to the image of a class 
$u\in H^{2n+2}(BG,\Z)$ by the map induced by the natural inclusion 
$\Z\subseteq\R$. We shall distinguish between two cases: In case I 
(the `Godbillon--Vey' case) we have $Q\in\ker(I^*(G)\to I^*(K))$, 
$K\subseteq G$ a maximal compact subgroup, and $u$ can be chosen to be $0$. 
Otherwise in case II we have $u\neq0$ (the `Cheeger--Chern--Simons case'). 
With this notation we shall prove the following in case I:

\begin{theorem}
  \label{1.2}
  Consider $Q\in I^{n+1}(G)$ as in case {\rm I} above and let $E\to Y$
  be a family of $G$-bundles with connections $\{A_z\mid z\in Z\}$ as
  in definition~$\ref{1.1}$. Let $\dim X=2n+1-\ell$ with
  $0\leq\ell\leq 2n+1$.

  \paritem{$(i)$} For $B$ a global extension of the family there is
  associated a natural class of $\ell$-forms
  $[\Lambda_{Y/Z}(Q,B)]\in\Omega^\ell(Z)/d\Omega^{\ell-1}(Z)$.
  
  \paritem{$(ii)$} This class is independent of the choice of extension
  provided $F^{n+1-\ell}_{A_z}=0$ for all $z\in Z$, where $F_{A_z}$ is
  the curvature form in the fibre $P_z$.

  \paritem{$(iii)$} Curvature formula $:$ ~ 
  \begin{align*}
    d\Lambda_{Y/Z}(Q,B)=(-1)^{\ell-1} \int_{Y/Z}Q(F_B^{n+1}), 
  \end{align*}
  where $Q(F_B^{n+1})\in\Omega^{2n+2}(Y)$ is the characteristic
  form associated to $Q$.
  
  \paritem{$(iv)$} If $F_{A_z}^{n-\ell}=0$ for all $z\in Z$ then
  $[\Lambda_{Y/Z}(Q,B)]$ lies in $H^\ell(Z,\R)$.
\end{theorem}

Here $\int_{Y/Z}$ denotes integration over the fibre in the bundle
$\pi\colon Y\to Z$. Also the \emph{curvature} $F_A$ of a connection
$A$ in a principal $G$-bundle $P\to X$ is defined as usual by
$F_A=dA+\tfrac12[A,A]$.

For $Q\in I_\Z^{n+1}(G)$ as in case II above we shall prove
(section~\ref{six}) a result analogous to Theorem~\ref{1.2} only the 
integral class $u\in H^{2n+2}(BG,\Z)$ has to be taken into account, 
and the deRham complex $\Omega^*(Z)$ is going to be replaced by the
\emph{simplicial} deRham complex (as in Dupont~\cite{D1} or~\cite{D2}) 
for the nerve of an open covering of $Z$. 
In terms of the above mentioned notion of gerbes with connections
(see section~\ref{two} below) we shall prove the following:

\begin{theorem}
  \label{1.3}
  Consider $Q\in I^{n+1}(G)$ and $u\in H^{2n+2}(BG,\Z)$ as in case
  {\rm II} above, and let $E\to Y$ be a family of $G$--bundles with
  connections $\{A_z\mid z\in Z\}$ as in definition~$\ref{1.1}$. Let
  $\dim X=2n+1-\ell$, $0\leq\ell\leq2n+1$.
  
  \paritem{$(i)$} For $B$ a global extension of the family there is
  associated a natural equivalence class of 
  $\ell$--gerbes $\theta=\theta(Q,u,B)$ with connection 
  $\omega=(\omega^0,\ldots,\omega^\ell)$ for a suitable open covering
  $\U=\{U_i\mid i\in I\}$.
  
  \paritem{$(ii)$} This class $[\theta,\omega]$ is independent of the
  choice of extension provided $F_{A_z}^{n+1-\ell}=0$ for all $z\in
  Z$, where $F_{A_z}$ is the curvature form in the fibre $P_z$.

  \paritem{$(iii)$} Curvature formula $:$ ~ 
  \begin{align}
    \label{1.4}
    d\omega^0=(-1)^{\ell-1}~\epsilon^*\int_{Y/Z}Q(F_B^{n+1})
    \qquad\hbox{and}\qquad
    \delta_*[\theta]=(-1)^{\ell-1}~\pi_!(u(E)).
  \end{align}

  \paritem{$(iv)$} If $F_{A_z}^{n-\ell}=0$ then $d\omega^0=0$ and the
  invariant $[\theta,\omega]$ lies in $H^\ell(Z,\R/\Z)$.
\end{theorem}

In \eqref{1.4} $\epsilon^*\colon\Omega^*(Z)\to\check C^0(\mathcal
U,\underline{\Omega}^*)$ is the natural inclusion of the deRham
complex into the \v Cech bicomplex. Furthermore $u(E)\in
H^{2n+2}(Y,\Z)$ is the associated characteristic class for the
$G$--bundle $E\to Y$ and $\pi_!\colon H^{2n+2}(Y,\Z)\to H^{\ell+1}(Z,\Z)$
is the usual transfer map. Finally
\begin{align*}
  \delta_*\colon H^\ell(Z,\underline{U(1)})\overset\cong \lra
  H^{\ell+1}(Z,\Z)
\end{align*}
is the usual isomorphism in \v Cech--cohomology. 

The above theorems contain the classical secondary characteristic
classes by taking $X=\{\pt\}$ and $\ell=2n+1$; but in this case
the invariants may depend on the extension $B$ (see section~\ref{five}). 
We are more concerned with the case $\ell\leq n$ where this does not happen. 
In particular we shall apply the Theorems~\ref{1.2} and~\ref{1.3} to 
families of foliated $G$--bundles of codimension $q$ in the sense of 
Kamber--Tondeur~\cite{KT}. These have \emph{adapted} connections $A$ whose
curvature $F_A$ satisfy $F_A^{q+1}=0$. Hence we obtain invariants for
families of such foliations provided $n-\ell\geq q$. We refer to
section~\ref{six} for a precise statement.

In the case $\ell=1$ Theorem~\ref{1.3} includes the construction of
the generalized Quillen line bundles considered in~\cite{DJ} which was
our motivating example. In section~\ref{six} we shall also consider a 
relative version of our construction generalizing the notion of a
`Chern-Simons section' considered in~\cite{DJ}.

Our Theorems~\ref{1.2} and~\ref{1.3} overlap with the results
of Freed~\cite{F3} but the methods are rather different. 
In fact we take advantage of the reformulation of 
`gerbes with connection' and smooth Deligne cohomology in terms
of simplicial differential forms as explained in section~\ref{three}.
In particular the notion of \emph{integration along the fibres} which
we are going to use, is fairly straight forward in this formulation
(see section 4 below or Dupont--Ljungmann~\cite{DL} ).
Also, as we shall see in section~\ref{five}, the Cheeger--Chern--Simons
characters are represented by simplicial differential forms. 
There are by now several ways of looking at gerbes with connection 
(see e.g. Hitchin~\cite{H}), but we hope to demonstrate that the
representation as a simplicial differential form is both an attractive 
and a convenient point of view.

The results of the paper go back a few years but the presentation
follows a talk given by the first author in November~2002 during the
program `Aspects of Foliation Theory' at the Erwin Schrödinger
Institute in Vienna. Both authors gratefully acknowledge the hospitality 
and support of the Erwin Schrödinger Institute. The second author visited 
{\AA}rhus on several occasions during the preparation of this work and 
would like to thank the Department of Mathematics at Aarhus University 
for its hospitality and support. 
Finally we want to thank the referee for some very useful comments in
particular on the terminology of `gerbes' and `Deligne cohomology'.


\medbreak
\section{Gerbes with connection}\label{two}

In this section we briefly recall the notion of a `gerbe with 
connection' and smooth `Deligne cohomology'. We refer to~\cite{B2}
for more information. We shall only consider \emph{Hermitian line 
gerbes} which are by definition \v Cech cocycles 
for the sheaf $\underline{U(1)}$ of smooth functions with values in 
the circle group $U(1)\subseteq\C$. 
For convenience we shall identify this group with
$\R/\Z$ via the map $z\leftrightarrow\frac1{2\pi i}\log z$, 
$z\in U(1)$. Hence a (Hermitian line) \emph{$p$--gerbe} on a smooth 
manifold $X$ is a $p$-cocycle in the \v Cech complex
\begin{align*}
  \check C^p(\U, \underline{\R/\Z})
  =\prod_{(i_0,\ldots,i_p)} C^\infty(U_{(i_0,\ldots,i_p)},\R/\Z),
\end{align*}
with the usual coboundary
\begin{align}
  \label{2.1}
  (\check\delta\theta)_{i_0,\ldots,i_p}
  =\sum_{\nu=0}^{p+1}(-1)^i\theta_{i_0,\ldots,\widehat{i_\nu},\ldots,i_p}.
\end{align}
Here $\U=\{U_i\mid i\in I\}$ is an open covering of $X$. For convenience 
we assume that $\U$ is `good' in the sense that all non-empty intersections 
$U_{i_0,\ldots,i_p}=U_{i_0}\cap\cdots\cap U_{i_p}$ 
are contractible. Is is well-known that every open covering has a good
refinement and that for such covering we have
\begin{align*}
  H^p(\check C^*(\U,\R/\Z))\cong
  H^p(X,\underline{\R/\Z}).
\end{align*}
Notice also that every cochain is the reduction of a cochain in
$\check C^*(\U,\underline\R)$ and that the isomorphism
\begin{align}
  \label{2.2}
  \delta_*\colon H^p(X,\underline{\R/\Z})\overset\cong\longrightarrow
  H^{p+1}(X,\Z)
\end{align}
is indeed induced by $\check\delta$ in \eqref{2.1} applied to such a
lift.

In general consider the \v Cech--deRham bicomplex
\begin{align}
  \label{2.3}
  \check\Omega_\R^{p,q}(\U)=\check C^p(\U,\underline{\Omega}^q)
\end{align}
with differential in the total complex $\check\Omega_\R^*(\U)$ given
on $\check\Omega_\R^{p,*}$ by $D=\check\delta+(-1)^pd$. Notice that
there are natural inclusions of chain complexes
\begin{align}
  \label{2.4}
  \check C^*(\U,\Z)\subseteq\check C^*(\U,\underline{\Omega}^0)
  \subseteq\check\Omega_\R^*(\U),
\end{align}
and
\begin{align}
  \label{2.5}
  \epsilon^*\colon\Omega^*(X)\overset\subseteq\longrightarrow
  \check C^0(\U,\underline{\Omega}^*)\subseteq
  \check\Omega_\R^*(\U),
\end{align}
where $\epsilon^*$ is induced by the natural map
\begin{align*}
  \epsilon\colon\bigsqcup_iU_i\to X.
\end{align*}
Since $\U$ is good we have
\begin{align*}
  \check C^*(\U,\underline{\R/\Z})=
  \check C^*(\U,\underline{\Omega}^0)/\check C^*(\U,\Z)
\end{align*}
and we put
\begin{align}
  \label{2.6}
  \check\Omega_{\R/\Z}^*(\U)=
  \check\Omega_\R^*(\U)/\check C^*(\U,\Z).
\end{align}
Notice that the canonical map 
\begin{align*}
  \epsilon^*\colon\Omega^*(X)\to\check\Omega_\R^*(\U)\to
  \check\Omega_{\R/\Z}^*(\U)
\end{align*}
is injective in degrees $>0$. We now have the following:

\begin{lemma}
  \label{2.7}
  Let $\U$ be a good covering of $X$. Then 

  \paritem{$(i)$} $H^*(\check\Omega_{\R}^*(\U)/\epsilon^*\Omega^*(X))=0$.
  
  \paritem{$(ii)$} $H^*(\check\Omega_{\R/\Z}^*(\U))\cong H^*(X,\R/\Z)$
  for $\R/\Z$ the constant sheaf.

  \paritem{$(iii)$} There is a natural isomorphism
  \begin{align*}
    D_*\colon
    H^\ell(\check\Omega_{\R/\Z}^*(\U)/\epsilon^*\Omega^*(X))\cong
    H^{\ell+1}(X,\Z)
  \end{align*}
  for $\ell\geq0$.
\end{lemma}

\begin{proof}
  \paritem{(i)} follows since
  $\epsilon^*\colon\Omega^*(X)\overset\subseteq\longrightarrow 
  \check\Omega_\R^*(\U)$ is a homology isomorphism. 

  \paritem{(ii)} follows since, for $\R$ the constant sheaf, the
  inclusion $\check C^*(\U,\R)\overset\subseteq\longrightarrow 
  \check\Omega_\R^*(\U)$ is a homology isomorphism. 

  \paritem{(iii)} Now $D_*$ is just the connecting homomorphism for
  the exact sequence
  \begin{align*}
    0\to
    \check C^*(\U,\Z)\to
    \check\Omega_{\R}^*(\U)/\epsilon^*\Omega^*(X)\to
    \check\Omega_{\R/\Z}^*(\U)/\epsilon^*\Omega^*(X)\to
    0.
  \end{align*}
\end{proof}

We can now define a gerbe with connection as follows:

\begin{definition}
  \label{2.8}
  Let $\U$ be a good covering for $X$.

  \paritem{(i)} A \emph{connection} $\omega$ in an $\ell$--gerbe 
  $\theta \in \check\Omega_{\R/\Z}^{\ell,0}(\U)~, ~\check\delta\theta = 0~,$ 
  is given by $\omega\in \check\Omega_{\R}^\ell(\U)$, that is a 
  sequence $\omega=(\omega^0,\ldots,\omega^\ell)~, ~
  \omega^\nu\in\check\Omega_{\R}^{\nu,\ell-\nu}(\U)$, 
  $\nu=0,\ldots,\ell$,
  with $\omega^\ell\equiv -\theta\mod\Z$, such that $\omega$ is a 
  cycle in $\check\Omega_{\R/\Z}^*(\U)/\epsilon^*\Omega^*(X)$.
  
  \paritem{(ii)} The \emph{curvature form} for $\omega$ is the unique
  closed $(\ell+1)$--form $F_\omega$ such that 
  $$
    \epsilon^* F_\omega=d\omega^0\in \check\Omega_{\R}^{0,\ell+1}(\U)
  $$
  The connection is called \emph{flat} if $F_\omega=0.$  
  
  \paritem{(iii)} Two $\ell$--gerbes $\theta_1,\theta_2$ with
  connections $\omega_1,\omega_2$ are \emph{equivalent} if
  $\omega_1-\omega_2$ is a coboundary in $\check\Omega_{\R/\Z}^*(\U)$.
  The set of equivalence classes $[\theta,\omega]$ is denoted
  $H^{\ell+1}_\D(X,\Z)$ and is called the (smooth) 
  \emph{Deligne cohomology} in degree $\ell+1$ (note the shift in degree).
\end{definition}

\begin{remarks}
\label{rem1}
  \paritem{1.} Thus $H^{\ell+1}_\D (X,\Z)$ is the homology of the sequence
  \begin{align}
    \label{2.9}
    \check\Omega_{\R/\Z}^{\ell-1}(\U)\overset d\lra
    \check\Omega_{\R/\Z}^\ell(\U)\overset d\lra
    \check\Omega_{\R/\Z}^{\ell+1}(\U)/\epsilon^*\Omega^{\ell+1}(X).
  \end{align}
  
  \paritem{2.} The set of equivalence classes of $\ell$--gerbes with
  \emph{flat} connections is isomorphic to $H^\ell(X,\R/\Z)$ by
  Lemma~\ref{2.7}.

  \paritem{3.} It follows also using Lemma~\ref{2.7}, that there is a
  natural exact sequence
  \begin{align}
    \label{2.11}
    0\to
    H^\ell(X,\R/\Z)\longrightarrow
    H_\D^{\ell+1}(X,\Z)\overset{d_*} \lra
    \Omega_{\rm cl}^{\ell+1}(X,\Z) \lra
    0.
  \end{align}
  Here $\Omega_{\rm cl}^{\ell+1}(X,\Z)\subseteq\Omega^*(X)$ denotes
  the set of closed forms with integral periods, and $d_*$ is induced
  by the map sending $\omega$ to the curvature form 
  $F_\omega$. 
  In particular, as the notation indicates,
  $H_\D^{\ell+1}(X,\Z)$ does not depend on the choice of a (good)
  covering $\U$.
  
  \paritem{4.} Notice the natural commutative diagram
  \begin{align}
    \label{2.12}
    \xymatrix{
      {H_\D^{\ell+1}(X,\Z)}
      \ar[r] \ar[d]&
      {H^\ell(X, \underline{\R/\Z})}
      \ar[d]_-{\cong}^-{\delta_*} \\
      {H^\ell(\check\Omega_{\R/\Z}^*(\U)/\epsilon^*\Omega^*(X))}
      \ar[r]^-{D_*}_-{\cong} &
      {H^{\ell+1}(X,\Z)}}
  \end{align}
  where the top horizontal map is induced by the map forgetting the
  connection and where $D_*$ is given by Lemma~\ref{2.7}.
  Note that by construction the invariants are defined by global 
  forms in case I, whereas the invariants are defined by simplicial 
  forms in case II. 
  
\paritem{5.} The explicit description of an $\ell$--gerbe $\theta$ 
  with connection $\omega$ is as follows. 
  Let $\omega$ be a sequence $(\omega^0,\ldots,\omega^\ell)$ of cochains 
  $\omega^\nu \in \check\Omega_{\R}^{\nu,\ell-\nu}(\U)~, ~
  \nu=0,\ldots,\ell$, satisfying
\begin{equation}\label{2.10}
    \begin{aligned}     
      \check\delta\omega^{\nu-1}&+(-1)^{\nu}d\omega^\nu = 0~,
      \qquad\nu=1,\ldots,\ell,\\      
      \check\delta\omega^\ell&\equiv0\mod\Z. 
    \end{aligned}
  \end{equation} 
  The first equation for $\nu=1$ in \eqref{2.10} implies that 
  $\check\delta d\omega^0=0$, and $d\omega^0$ defines a global closed 
  $(\ell+1)$--form $F_\omega$, that is $\epsilon^* F_\omega=
  d\omega^0$. 
  The last equation 
  in \eqref{2.10} says that $-\check\delta\omega^\ell$ is an integral 
  $(\ell+1)$--cycle $z \in \check Z^{\ell+1}(\U,\Z)$, that is 
  $-\omega^{\ell} \in \check\Omega_{\R}^{\ell,0}(\U)$ is the lift 
  of a unique $\ell$--cycle $\theta \in \check\Omega_{\R/\Z}^{\ell,0}(\U)$. 
  Thus from \eqref{2.2} we have $\delta_* [\theta] = [z]$. 
  Moreover by construction, the integral class $[z]$ determines the class 
  $[F_\omega]$ under the canonical homomorphism 
  $r:H^{\ell+1}(X,\Z) \to H^{\ell+1}(X,\R)$. 
   Then $\omega$ is a connection for the $\ell$--gerbe $\theta$. 

\paritem{6.} In terms of the notation in~\cite{B2}
our smooth Deligne cohomology group $H^{\ell+1}_\D(X,\Z)$ is canonically 
isomorphic to the group $H^{\ell+1}_{\D,\infty} (X,\Z(\ell+1))$, that is 
the hypercohomology group of $X$ in degree $\ell+1$ 
with values in the sheaf complex 
$$
\Z \ra \underline{\Omega}^0 \ra \underline{\Omega}^1 \ra \ldots 
\ra \underline{\Omega}^\ell~. 
$$ 
Since in the smooth case $H^k_{\D,\infty} (X,\Z(\ell+1))$
is ordinary cohomology with coefficients $\R / \Z$ for $k<{\ell+1}$, 
respectively $\Z$ for $k>{\ell+1}$, 
$k={\ell+1}$ is the only degree which needs a special name 
and we have therefore deleted the extra index from the notation. 
This is of course in contrast to the holomorphic Deligne cohomology for 
an algebraic variety.
\end{remarks}

Finally let us mention the interpretation of $H_\D^*(X,\Z)$ as the
group of \emph{differential characters} in the sense of
Cheeger-Simons~\cite{CrS} (see also Dupont et.al.~\cite{DHZ}). 
Let $C_*^\Sing(X)$ denote the chain complex of (smooth) singular 
chains in $X$ and let 
\begin{align*}
  \I\colon\Omega^*(X)\to C_\Sing^*(X,\R)=\Hom_\Z(C_*^\Sing(X),\R)
\end{align*}
be the deRham integration map.

\begin{definition}
  \label{2.13}
  The group of \emph{differential characters} (mod $\Z$) in degree
  $\ell+1$ is
  \begin{align*}
    \what H^{\ell+1}(X,\Z)
    =\{(f,\alpha)\in\Hom_\Z(Z_\ell^\Sing(X),\R/\Z)\oplus
    \Omega^{\ell+1}(X)\mid\delta
    \hbox{$f=\I(\alpha)$ and $d\alpha=0$}\}.
  \end{align*}
\end{definition}
Here $Z_\ell^\Sing(X)\subseteq C_\ell^\Sing(X)$ is the set of
cycles. The following is well-known (cf.~\cite{DHZ}) but is included for
completeness:

\begin{proposition}
  \label{2.14}
  There is a natural isomorphism $H_\D^*(X,\Z)\cong\what H^*(X,\Z)$.
\end{proposition}

\begin{proof}
  Choose a good open covering $\U=\{U_i\mid i\in I\}$ of $X$ and let
  $i\colon C_*^\Sing(X,\U)\subseteq C_*^\Sing(X)$ be the inclusion of
  the subcomplex generated by $\bigcup_{i\in I}C_*^\Sing(U_i)$. Since
  $i$ is a chain equivalence we can choose a chain map
  \begin{align*}
    p\colon C_*^\Sing(X)\to C_*^\Sing(X,\U)
  \end{align*}
  such that $p\circ i=\id$ and $i\circ p$ is chain homotopic to the
  identity with chain homotopy $s$. Then for $(f,\alpha)\in
  \what H^{\ell+1}(X,\Z)$ and $\xi\in Z_\ell^\Sing(X)$ we have
  \begin{align*}
    \langle f,\xi\rangle-\langle f,i\circ p(\xi)\rangle
    =\langle\delta f,s(\xi)\rangle
    =\langle \I(\alpha),s(\xi)\rangle
  \end{align*}
  so that $f$ is determined by its restriction to the set of cycles
  $Z_\ell^\Sing(X,\U)$ in the chain complex $C_*^\Sing(X,\U)$. Hence
  we can replace $Z_\ell^\Sing(X)$ by $Z_\ell^\Sing(X,\U)$ in
  Definition~\ref{2.13}. Now we consider the \v Cech bicomplex of
  singular chains
  \begin{align*}
    \check C_{p,q}^\Sing(\U)=\bigoplus_{(i_0,\ldots,i_p)}C_q^\Sing 
    (U_{i_0 \cdots i_p})
  \end{align*}
  with associated total complex $\check C_*^\Sing(\U)$. Then again the
  natural chain map
  \begin{align*}
    \xymatrix@C-3em{
      {\check C_*^\Sing(\U)}
      \ar[rr]^-{\epsilon_*}\ar[dr]&&
      {C_*^\Sing(X,\U)}\\
      &{\check C_{0,*}^\Sing(\U)}
      \ar[ur]}
  \end{align*}
  induced by $\epsilon\colon\bigsqcup_{i\in I}U_i\to X$, has an
  `inverse' chain map $j$ such that $\epsilon_*\circ j=\id$ and
  $j\circ\epsilon_*$ is chain homotopic to the identity. Now we can
  define a map 
  \begin{align*}
    j_*\colon H_\D^*(X,\Z)\to\what H^*(X,\Z)
  \end{align*}
  by $j_*[\omega,\theta]=(f,\alpha)$ where $f(\xi)=\langle
  \I(\omega),j(\xi)\rangle$, $\xi\in Z_\ell^\Sing(X,\U)$, and
  $\alpha=(\epsilon^*)\inv d\omega^0$. In fact for $x\in
  C_{\ell+1}^\Sing(X,\U)$ we have
  \begin{align*}
    \langle\delta f,x\rangle
    &=\langle \I(\omega),\partial j(x)\rangle
    =\langle \I(D(\omega)),j(x)\rangle\\
    &=\langle \I(d\omega^0),j(x)\rangle
    =\langle \I(\alpha),\epsilon_*j_*(x)\rangle
    =\langle \I(\alpha),x\rangle
  \end{align*}
  so that $(f,\alpha)\in\what H^{\ell+1}(X,\Z)$. Since any two choices
  of $j$ are chain homotopic, it is also straight forward to see that
  $j_*$ does not depend on the particular choice. Finally, in order to
  show that $j_*$ is an isomorphism one just observes that there is a
  natural exact sequence similar to the one in \eqref{2.11}:
  \begin{align}
    \label{2.15}
    0\to
    H^\ell(X,\R/\Z)\to
    \what H^{\ell+1}(X,\Z)\to
    \Omega_{\rm cl}^{\ell+1}(X,\Z)\to
    0
  \end{align}
  where the second map is the one sending $(f,\alpha)$ to $\alpha$.
\end{proof}


\medbreak
\section{Gerbes and simplicial forms}\label{three}

In this section we shall reformulate the smooth Deligne cohomology in
terms of simplicial deRham cohomology as in \cite{D1}, \cite{D2} and 
\cite{DJt}. 
As before let $X$ be a smooth manifold and let 
$\U=\{U_i\}_{i\in I}$ be a good covering of $X$. 
For convenience we choose a linear
ordering of the index set $I$. The \emph{nerve} $N\U$ of $\U$ is the
simplicial manifold $N\U=\{N\U(p)\}_{p\geq0}$, given by 
\begin{align}
  \label{3.1}
  N\U(p)=\coprod_{i_0\leq\cdots\leq i_p}U_{i_0\cdots i_p}
\end{align}
and with face and degeneracy operators $\epsilon_i\colon N\U(p)\to
N\U(p-1)$, $i=0,\ldots,p$, $\eta_j\colon N\U(p)\to N\U(p+1)$,
$j=1,\ldots,p$, given by the obvious inclusion maps corresponding to
deletion of the $i$-th index respectively repeating the $j$-th
index. Also let $\Delta^p\subseteq\R^{p+1}$ be the \emph{standard}
$p$-simplex
\begin{align*}
  \Delta^p=\Big\{t=(t_0,\ldots,t_p)\Bigm| t_i\geq0,\sum_{i=0}^pt_i=1\Big\}
\end{align*}
with the corresponding face and degeneracy maps
$\epsilon^i\colon\Delta^{p-1}\to\Delta^p$, $i=0,\ldots,p$,
respectively $\eta^j\colon\Delta^{p+1}\to\Delta^p$, $j=0,\ldots,p$.

\begin{definition}
  \label{3.2}
  \paritem{(i)} A \emph{simplicial} $k$-form $\omega$ on $N\U$ is a
  sequence of $k$-forms $\omega^{(p)}$ on $\Delta^p\times N\U$
  satisfying 
  \begin{align*}
    (\epsilon^i\times\id)^*\omega^{(p)}=(\id\times\epsilon_i)^*\omega^{(p-1)},
    \qquad i=0,\ldots,p,\quad p=1,2,\ldots.
  \end{align*}
  
  \paritem{(ii)} $\omega$ is called \emph{normal} if it furthermore
  satisfies
  \begin{align*}
    (\eta^i\times\id)^*\omega^{(p)}=(\id\times\eta_i)\omega^{(p+1)},
    \qquad i=1,\ldots,p,\quad p=1,2,\ldots.
  \end{align*}
\end{definition}

We shall denote the set of simplicial $k$-forms (respectively normal
$k$-forms) by $\Omega^k(||N\U||)$ (respectively $\Omega^k(|N\U|)$)
corresponding to the `fat' (respectively `thin') realizations
$||N\U||$ (respectively $|N\U|$). Clearly $\Omega^*(||N\U||)$ 
is a differential graded algebra and
$\Omega^*(|N\U|)\subseteq\Omega^*(||N\U||)$ is a
DGA-subalgebra. Notice that the inclusions $U_i\subseteq X$ induce a
natural simplicial map $\epsilon\colon N\U\to N\{X\}$ and this in turn
induces a DGA-map
\begin{align}
  \label{3.3}
  \epsilon^*\colon\Omega^*(X)\to\Omega^*(|N\U|)\subseteq\Omega^*(||N\U||)
\end{align}
where $\Omega^*(X)=\Omega^*(|N\{X\}|)$ is the usual deRham complex. It
follows from \cite{D1} that $\epsilon^*$ induces homology isomorphisms
\begin{align}
  \label{3.4}
  \epsilon^*\colon H(\Omega^*(X))\underset\cong\longrightarrow
  H(\Omega^*(|N\U|))\underset\cong\longrightarrow
  H(\Omega^*(||N\U||)).
\end{align}

The relation with the \v Cech--deRham complex in section~\ref{two} 
is given by the integration map 
\begin{align}
  \label{3.5}
  \I_\Delta\colon\Omega^{p,q}(||N\U||)\to\check\Omega_\R^{p,q}(\U),
  \quad \I_\Delta(\omega)=\int_{\Delta^p}\omega^{(p)},
\end{align}
where $\omega$ lies in $\Omega^{p,q}$ if it has degree $p$ as a form
in the variables of $\Delta^n$, $n\geq p$. This is a map of
bicomplexes and again by \cite{D1} the corresponding map of 
total complexes induces an isomorphism
\begin{align}
  \label{3.6}
  \I_\Delta\colon H(\Omega^*(||N\U||))\underset\cong\longrightarrow
  H(\check\Omega_\R^*(\U)).
\end{align}
Also $\I_\Delta$ clearly commutes with $\epsilon^*$ given by
\eqref{2.5} and \eqref{3.3}. 

For the representation of the integral cohomology we also consider the
discrete simplicial set $N_d\U$ where a $p$-simplex is a point
$(i_0,\ldots,i_p)$ for each non-empty intersection
$U_{i_0}\cap\cdots\cap U_{i_p}$, $i_0\leq i_1\leq\ldots\leq i_p$, and
we let $\eta\colon N\U\to N_d\U$ denote the simplicial map sending
$U_{i_0}\cap\cdots\cap U_{i_p}$ to $(i_0,\ldots,i_p)$. Notice that for
$\U$ a good covering we have a commutative diagram of homotopy
equvalences
\begin{align}
  \label{3.7}
  \vcenter{\xymatrix{
    {||N\U||}\ar[r]_-\simeq\ar[d]_-{||\eta||}&
    {|N\U|}\ar[d]^-{|\eta|}\\
    {||N_d\U||}\ar[r]_-\simeq&
    {|N_d\U|}}}
\end{align}
and a similar diagram of isomorphisms
\begin{align}
  \label{3.8}
  \vcenter{\xymatrix{
    {H(\Omega^*(||N\U||))}&
    {H(\Omega^*(|N\U|))}\ar[l]^-\cong\\
    {H(\Omega^*(||N_d\U||))}\ar[u]_-\cong^-{\eta^*}&
    {H(\Omega^*(|N_d\U|))}\ar[l]^-\cong\ar[u]_-{\eta^*}^-\cong}}
\end{align}

Also notice that $\eta^*$ maps 
\begin{align*}
  \Omega^*(||N_d\U||)=\Omega^{*,0}(||N_d\U||)
\end{align*}
injectively into $\Omega^{*,0}(||N\U||)\subseteq\Omega^*(||N\U||)$ and
that $\omega\in\Omega^*(||N\U||)$ lies in the image if and only if it
only involves the variables of $\Delta^p$.

\begin{definition}
  \label{3.9}
  
  \paritem{(i)} A $k$-form $\omega\in\Omega^*(||N\U||)$ is called
  \emph{discrete} if $\omega\in\eta^*(\Omega^*(||N_d\U||)$.

  \paritem{(ii)} $\omega\in\Omega^*(||N\U||)$ is called
  \emph{integral} if it is discrete and if furthermore 
  \begin{align*}
    \I_\Delta(\omega)\in\check
    C^*(\U,\Z)\subseteq\check\Omega^{*,0}(\U).
  \end{align*}
\end{definition}

We let $\Omega^*_\Z(||N\U||)\subseteq\Omega^*(||N\U||)$ (respectively
$\Omega^*_\Z(|N\U|)\subseteq\Omega^*(|N\U|)$) denote the chain complex
of integral forms (respectively integral normal forms) and we also put
\begin{align}
  \label{3.10}
  \Omega^*_{\R/\Z}(||N\U||)=\Omega^*(||N\U||)/\Omega^*_\Z(||N\U||)
\end{align}
respectively
\begin{align}
  \label{3.11}
  \Omega^*_{\R/\Z}(|N\U|)=\Omega^*(|N\U|)/\Omega^*_\Z(|N\U|).
\end{align}
We now have the following:

\begin{proposition}
  \label{3.12}
  Let $\U$ be a good covering. Then there are natural isomorphisms

  \paritem{$(i)$} $H(\Omega^*_\Z(||N_d\U||))\overset{\eta^*}\cong
  H(\Omega^*_\Z(||N\U||))\overset{{\I_\Delta}_*}\cong H(\check
  C^*(\U,\Z))=H^*(X,\Z)$,
  
  \paritem{$(ii)$}
  $H(\Omega^*_{\R/\Z}(||N\U||))\overset{{\I_\Delta}_*}\cong H(\check
  \Omega^*_{\R/\Z}(\U))\cong H^*(X,\R/\Z)$,

  \paritem{$(iii)$}
  $H^\ell(\Omega^*(||N\U||)/(\Omega^*_\Z(||N\U||)+\epsilon^*\Omega^*(X)))
  \overset{d_*}\cong H^{\ell+1}(\Omega^*_\Z(||N\U||))\cong
  H^{\ell+1}(X,\Z)$.

  \paritem{$(iv)$} Furthermore $\I_\Delta$ induces a natural isomorphism
  to $H_\D^{\ell+1}(X,\Z)$ from the homo\-logy of the sequence
  \begin{align}
    \label{3.13}
    \Omega^{\ell-1}_{\R/\Z}(||N\U||)\overset d\longrightarrow
    \Omega^\ell_{\R/\Z}(||N\U||)\overset d\longrightarrow
    \Omega^{\ell+1}_{\R/\Z}(||N\U||)/\epsilon^*\Omega^{\ell+1}(X).
  \end{align}

  \paritem{$(v)$}
  In $(i)$--$(iv)$ above $||N\U||$ can be replaced by $|N\U|$. 
\end{proposition}

\begin{proof}
  \paritem{(i)} In the commutative diagram
  \begin{align*}
    \xymatrix@R-2em{
      {\Omega^*_\Z(||N_d\U||)}
      \ar[dr]^-{\I_\Delta}\ar[dd]_-{\eta^*}\\
      &{\check C^*(\U,\Z)}\\
      {\Omega^*_\Z(||N\U||)}
      \ar[ur]_-{\I_\Delta}}
  \end{align*}
  $\eta^*$ is an isomorphism and $\I_\Delta$ for $N_d\U$ is a homology
  isomorphism since it is surjective and the kernel has vanishing
  homology by the simplicial deRham theorem. Hence also $\I_\Delta$
  for $N\U$ is a homology isomorphism.

  \paritem{(ii)} now follows from (i) and \eqref{3.6} together with
  Lemma~\ref{2.7} (ii).
  
  \paritem{(iii)} is similar to lemma~\ref{2.7}, (iii). 

  \paritem{(iv)} follows from the five--lemma applied to the sequence
  in \eqref{2.11} and the corresponding sequence for the homology
  group in \eqref{3.13}.

  \paritem{(v)} follows similarly. 
\end{proof}

\begin{corollary}
  \label{3.13a} 
  Every class in
  $H_\D^{\ell+1} (X,\Z)$ can be represented by an $\ell$--gerbe $\theta$
  with connection $\omega$ of the form $\omega=\I_\Delta(\Lambda)$ for
  some simplicial $\ell$-form $\Lambda\in\Omega^\ell(||N\U||)$
  satisfying
  \begin{align}
    \label{3.14}
    d\Lambda=\epsilon^*\alpha-\eta^*\beta,
    \qquad\alpha\in\Omega^{\ell+1}(X),
    \quad\beta\in\Omega^{\ell+1}_\Z(||N_d\U||).
  \end{align}
  Furthermore $\Lambda$ and $\beta$ can be chosen to be normal in the
  sense of Definition~$\ref{3.2}$. 
\end{corollary}

\begin{remarks}
\label{rem2}
  \paritem{1.} We shall call a (normal) simplicial
  $\ell$-form $\Lambda$ a (\emph{normal}) \emph{simplicial
    $\ell$--gerbe} if it satisfies \eqref{3.14}.
  
  \paritem{2.} Continuing with the previous notation, we write
  \begin{align*}
    \Lambda=\Lambda^0+\cdots+\Lambda^\ell\in
    \bigoplus_{\nu=0}^\ell\Omega^{\nu,\ell-\nu}(||N\U||)
  \end{align*}
  and we put
  \begin{align}
    \label{3.15}
    \theta= -\int_{\Delta^\ell}\Lambda^\ell,\qquad
    \omega^\nu=\int_{\Delta^\nu}\Lambda^\nu,\qquad
    \nu=0,\ldots,\ell.
  \end{align}
  Then \eqref{3.14} corresponds to the condition \eqref{2.10} for the
  $\ell$--gerbe $\theta$ with connection
  $\omega=(\omega^0,\ldots,\omega^\ell \equiv -\theta)$.

  \paritem{3.} Note that $\alpha$ and $\beta$ in \eqref{3.14} are
  uniquely determined by $\Lambda$ and that $\alpha$ is the curvature
  form of $\omega$. We shall refer to it as the \emph{curvature form}
  for $\Lambda$.

  \paritem{4.} By \eqref{3.14} and \eqref{3.15} we have
  \begin{align}
    \label{3.16}
    \I_\Delta(\beta)=-\int_{\Delta^{\ell+1}}d\Lambda^\ell
    =\check\delta\theta\in\check C^{\ell+1}(\U,\Z).
  \end{align}
  Hence $\beta$ represents the characteristic class
  $$
  z=\check\delta_*[\theta] \in
  H^{\ell+1}(X,\Z)=H^{\ell+1}(\Omega^*_\Z(||N_d\U||)).
  $$

  \paritem{5.} The simplicial deRham complexes $\Omega^*(||N\U||)$
  and $\Omega^*(|N\U|)$ as well as the corresponding subcomplexes of
  integral forms are clearly functorial with respect to smooth maps
  $f\colon X'\to X$ and \emph{compatible coverings}. By this we mean
  coverings $\U'=\{U'_{i'}\}_{i'\in I'}$ of $X'$ and $\U=\{U_i\}_{i\in
  I}$ of $X$ together with an order preserving map $\nu\colon I'\to I$
  such that $f(U'_{i'})\subseteq U_{\nu(i')}$ for all $i'\in I'$;
  that is, $\U'$ is a refinement of $f\inv(\U)$. The induced maps in the
  deRham complexes do depend on $\nu$ but the induced map in Deligne
  cohomology does not. Notice that this is the case also for
  $f=\id\colon X\to X$, that is, when $\U'$ is a refinement of $\U$. 

  \paritem{6.} If $X$ has dimension $m$ then it also has covering
  dimension $m$ (see e.g. \cite{HW}, chap.~II). Hence by taking a suitable
  refinement we obtain a covering $\U'$ for which $N\U'$ has only
  \emph{non-degenerate} simplices of dimension $\leq m$. In particular
  for such a covering we have 
  \begin{align}
    \label{3.17}
    \Omega^{k,\ell}(|N\U'|)=0
    \qquad\hbox{and}\qquad
    \Omega^k(|N_d\U'|)=0
    \qquad\hbox{for $k>m$}.
  \end{align}
\end{remarks}


\medbreak
\section{Fibre integration of simplicial forms}\label{four}

Fibre integration in smooth Deligne cohomology can be done in various
ways, see e.g. Freed~\cite{F3}, Gomi--Terashima~\cite{GT} or 
Hopkins--Singer~\cite{HS}.
In this section we sketch how to define it in terms of simplicial
forms. We refer to Dupont--Ljungmann~\cite{DL} for the details.

In the following $X$ denotes an oriented compact manifold of dimension
$m$ possibly with boundary and $\pi\colon Y\to Z$ is a smooth fibre
bundle with fibre $X$ and structure group $\Diff^+(X)$ of orientation
perserving diffeomorphisms. Also let $\V=\{V_j\}_{j\in J}$ and
$\U=\{U_i\}_{i\in I}$ be good open coverings of $Y$ respectively $Z$
(not necessarily compatible). We shall define \emph{integration along
  the fibre} for a \emph{normal} simplical $(k+m)$-form
$\omega\in\Omega^{k+m}(|N\V|)$ as a simplicial $k$-form
$\int_{Y/Z}\omega\in\Omega^k(||N\U||)$ defined by usual fibre
integration in the bundle $\Delta^p\times
N(\pi\inv\U)(p)\to\Delta^p\times N\U(p)$, $p=0,1,2,\ldots$ with fibre
$X$:
\begin{align}
  \setbox0\hbox{%
  $\scriptstyle(\Delta^p\times N(\pi\inv\U)(p))/(\Delta^p\times N\U(p))$}
  \label{4.1}
  \int_{Y/Z}\omega|_{\Delta^p\times N(\pi\inv\U)(p)}
  =\hskip.5\wd0\int\limits_{\hbox to0pt{\hss\box0\hss}}
  \tilde\phi^*\omega,
\end{align}
where $\pi\inv\U=\{\pi\inv U_i\}_{i\in I}$ is the obvious covering of
$Y$ and $\tilde\phi\colon||N(\pi\inv\U)||\to|N\V|$ denotes a
`piecewise smooth' map associated to a choice of partition of unity
for the coverings $\{\pi\inv U_i\cap V_j\}_{j\in J}$ for each $i\in
I$. For the construction of $\tilde\phi$ let us assume for simplicity
that $\pi\colon Y\to Z$ is the product fibration $X\times Z\to Z$. For
the case of a general fibration we refer to \cite{DL}. 
By remark~\ref{rem2}, 5 we can assume that 
$\V=\U'\times\U=\{V_{ij}=U'_j\times U_i\}_{i\in I,\;j\in J}$ where
$\U=\{U_i\}_{i\in I}$ and $\U'=\{U_j'\}_{j\in J}$ are open coverings
af $Z$ and $X$ respectively and we order $I\times J$ lexicographically
with $i\in I$ before $j\in J$. (Notice the interchange in $V_{ij}$.)
Also as in remark~\ref{rem2}, 6 we can assume that $N\U'$ has only 
non-degenerate simplices of dimension $\leq m$ and that 
$N(\U'\cap\partial X)$ has only non-degenerate simplices of dimension 
$\leq m-1$ ($m=$ dimension of $X$). Finally we choose a partition of 
unity $\{\phi_j\}_{j\in J}$ subordinate $\U'$. 
Then the natural projection $|N\U'|\to X$ has a
right-inverse $\bar\phi\colon X\to|N\U'|$ defined by
\begin{align}
  \label{4.2}
  \bar\phi(x)=((\phi_{j_0}(x),\ldots,\phi_{j_q}(x),x)_{j_0\cdots j_q}
  \in\Delta^q\times N\U'(q)
\end{align}
for those $x\in U_{j_0\cdots j_q}\subseteq X$ satisfying
$\phi_{j_0}(x)+\cdots+\phi_{j_q}(x)=1$. Now, we would like to define
$\tilde\phi$ in a similar fashion as the composite in the diagram
\begin{align}
  \label{4.3}
  \vcenter{
    \xymatrix@C+1em{
      {||N(\pi\inv\U)||}\ar[dr]^-{\tilde\phi}\ar[d]\\
      {|N(\pi\inv\U)|}\ar@{=}[d]&
      {|N\U'\times N\U|\hbox to0pt{${}=|N(\U'\times\U)|$\hss}}
      \ar[d]^-\approx_-\tau\\
      {X\times|N\U|}\ar[r]^-{\bar\phi\times\id}&
      {|N\U'|\times|N\U|}}}
\end{align}
where the homeomorphism $\tau$ is induced by the Eilenberg-Zilber
triangulation map
\begin{align*}
  \Delta^n\times(N\U'(n)\times N\U(n))\to
  (\Delta^n\times N\U'(n))\times(\Delta^n\times N\U(n))
\end{align*}
given by the diagonal $\Delta^n\to\Delta^n\times\Delta^n$.
It is well-known that $\tau\inv$ is given by the triangulation of a
prism $\Delta^q\times\Delta^p$ into $n$-simplices ($n=p+q$) one for
each `$(q,p)$-shuffle' of $(0,\ldots,n)$, that is, a pair of
non-decreasing functions 
\begin{align*}
  (\nu,\mu)\colon\{0,\ldots,n\}\to\{0,\ldots,q\}\times\{0,\ldots,p\}
\end{align*}
satisfying
\begin{align}
  \mu(0)=\nu(0)=0,\qquad\mu(n)=p,\qquad\nu(n)=q,\qquad\hbox{and}\label{4.4}\\
  \mu(r)-\mu(r-1)+\nu(r)-\nu(r-1)=1,\qquad r=1,\ldots,n,\label{4.5}
\end{align}
(so that for increasing $r$ the functions $\mu$ and $\nu$ alternate
increasing by $1$). It follows that
$\tilde\phi^*\omega\in\Omega^{k+m}(||N(\pi\inv\U)||)$ is the
simplicial form defined explicitly on $\Delta^p\times(X\times
U_{i_0\cdots i_p})$ in a neighborhood of a point $(t,x,z)$ by the sum
\begin{align}
  \label{4.6}
  (\tilde\phi^*\omega)_{i_0\cdots i_p}
  =\sum_{(\nu,\mu)}\tilde\phi^*_{(\nu,\mu)}\omega
\end{align}
with $(\nu,\mu)$ running through the $(q,p)$-shuffles as above. Here
$q$ is determined such that $\phi_{j_0}+\cdots+\phi_{j_q}=1$ near $x$
and
\begin{align*}
  \tilde\phi_{(\nu,\mu)}\colon
  \Delta^p\times(U'_{j_0\cdots j_q}\times U_{i_0\cdots,i_p})\to
  \Delta^n\times(U'_{j_{\nu(0)}}\times U_{i_{\mu(0)}})\cap\cdots\cap
  (U'_{j_{\nu(n)}}\times U_{i_{\mu(n)}})
\end{align*}
is given by the formula
\begin{align}
  \label{4.7}
  \tilde\phi_{(\nu,\mu)}(t,x,z)=(\sigma_0,\ldots,\sigma_n,x,z)
\end{align}
where
\begin{align}
  \label{4.8}
  \sigma_r=\sum_{(\nu',\mu')}t_{\mu'}\phi_{j_{\nu'}}(x)
\end{align}
is a sum over the pairs of integers $(\nu',\mu')$, $\mu'=1,\ldots,p$,
$\nu'=0,\ldots,q$, satisfying
$(\nu(r-1),\mu(r-1))<(\nu',\mu')\leq(\nu(r),\mu(r))$ in the
lexicographical order.  That is,
\begin{align*}
    \sigma_r=
    \begin{cases}
    t_{\mu(r)}\phi_{j_{\nu(r)}}(x) 
    &\hbox{if $\mu(r-1)=\mu(r)$, $\nu(r-1)<\nu(r)$},\\
    \noalign{\bigskip}
    { \begin{aligned} 
    & t_{\mu(r-1)}\sum_{\nu(r)<\nu'}\phi_{j_{\nu'}}(x)~+~ \\
     \noalign{\medskip} 
    &+~ t_{\mu(r)}\sum_{\nu'\leq\nu(r)}\phi_{j_{\nu'}}(x) 
     \end{aligned} }
    &\hbox{if $\mu(r-1)<\mu(r)$, $\nu(r-1)=\nu(r)$}.
   \end{cases}
\end{align*}

The form given by \eqref{4.6} clearly defines a smooth form in
$\Delta^p\times N(\pi\inv\U)(p)$ so that $\int_{Y/Z}\omega$ is indeed
well-defined by the formula \eqref{4.1}. Also it is easy to see from
the construction that it is a simplicial $k$-form, i.e., that it
satisfies Definition~\ref{3.1} (i). It is however not necessarily a
\emph{normal} simplicial form even though $\omega$ was normal to begin
with.

We note the following properties of fibre integration. The signs are
determined by the convention that we always integrate the variables
starting from the left:

\begin{proposition}
  \label{4.10}
  \paritem{$(i)$} Let $\omega\in\Omega^{k+m-1}(|N\V|)$, $m=\dim X$. Then
  \begin{align*}
    \int_{Y/Z}d\omega=\int_{\partial Y/Z}\omega
    +(-1)^{m}~d\int_{Y/Z}\omega.
  \end{align*}
  
  \paritem{$(ii)$} If $\partial X=\emptyset$ and
  $\omega\in\Omega_\Z^*(|N\V|)$ then $\int_{Y/Z}\omega$ is also
  integral.

  \paritem{$(iii)$} Suppose $\partial X=\emptyset$. Then
  \begin{align*}
    \int_{Y/Z}\colon\Omega^{k+m}(|N\V|)\to\Omega^k(||N\U||)
  \end{align*}
  induces the usual transfer map $\pi_!\colon H^{k+m}(Y)\to H^k(Z)$
  with coefficients in $\R$, $\Z$ or $\R/\Z$. Also it induces a
  well-defined map of smooth Deligne cohomology
  \begin{align*}
    \pi_!\colon H_\D^{k+m}(Y,\Z)\to H_\D^k(Z,\Z)
  \end{align*}
  independent of choices of coverings and partition of unity.

  \paritem{$(iv)$} $\pi_!$ is functorial with respect to bundle maps and
  compatible coverings.
\end{proposition}

\begin{proof}
  Again we restrict to the product case, referring to \cite{DL} for the
  general case.

  \paritem{(i)} By \eqref{4.1} this follows as for usual fibre
  integration from Stokes' Theorem.
  
  \paritem{(ii)} We shall prove that the \v Cech cochain
  $c=\I_\Delta\int_{Y/Z}\omega$ for the covering $\U$ has integral
  values. For this we observe that $(\tilde\phi^*\omega)_{i_0\cdots
    i_k}$ in \eqref{4.6} only involves $\omega$ restricted to the
  $(k+m)$-skeleton of $N_d\V$, hence by \eqref{3.17} can be assumed to
  be a closed integral form. Since $\bar\phi\colon X\to|N\U'|$ has
  degree one, it is straight forward from \eqref{4.3} and the
  Eilenberg--Zilber Theorem that $c_{i_0\cdots i_k}$ is the evaluation
  of $\I_\Delta(\omega)$ on the chain $[X]\times(i_0,\ldots,i_p)$ and
  hence is integral. In fact it follows that $c$ represents the slant
  product $\I_\Delta(\omega)/[X]$ in the integral cohomology.

  \paritem{(iii)} Since $\pi_!\colon H^{k+m}(Y)\to H^k(Z)$ is induced
  by the slant product by $[X]$ the first statement is already
  contained in the proof of (ii). That $\pi_!$ in Deligne cohomology
  is independent of choice of partition of unity follows from (i)
  applied to $Y\times[0,1]$ and the partition of unity
  $\{(1-t)\phi_j+t\phi'_j\}_{j\in J}$ where $\{\phi_j\}_{j\in J}$ and
  $\{\phi'_j\}_{j\in J}$ are the two given ones for the covering
  $\V$. Independence of choice of covering is now straightforward
  using remark~\ref{rem2}, 5. 

  \paritem{(iv)} is also straightforward.
\end{proof}

\begin{remark}
\label{rem3}
  In the case of a product fibration $\pi\colon X\times Z\to Z$ with
  the covering $\V=\{U'_j\times U_i\}_{(i,j)\in I\times J}$ as above, 
  we can also 
  represent a class in $H_\D^*(X\times Z,\Z)$ by
  a normal simplicial form $\omega$ in the bisimplicial manifold
  $N\U'\times N\U$ (cf. \cite{DJt}), i.e., by a collection of compatible
  forms on $\Delta^q\times\Delta^p\times N\U'(q)\times N\U(p)$. We can
  then define $\int_\xi\omega\in\Omega^{n-\ell}(|N\U|)$ for
  $\xi\in\check C_\ell^\Sing(N\U')$ any class in the \v Cech bicomplex
  of singular chains in the notation at the end of section~\ref{two} above. 
  In fact, for a singular $r$-simplex 
  $\xi=\sigma\colon\Delta^r\to U'_{j_0\cdots j_q}\subseteq X$ with $r+q=\ell$,  
  we just integrate the pull-back of $\omega$ to 
  $\Delta^q\times\Delta^p\times\Delta^r\times N\U(p)$ over 
  $\Delta^q\times\Delta^r$. For $\xi=[X]$ a representative for the 
  fundamental cycle of $X$ (in case $\partial X=\emptyset$) we have 
  (with $\tau$ being the Eilenberg-Zilber map as in \eqref{4.3} above):
  \begin{align}
    \label{4.11}
    \int_{[X]}\omega=\int_{Y/Z}\tau^*\omega.
  \end{align}
  Also we have the Stokes' formula similar to Proposition~\ref{4.10}~(i):
  \begin{align}
    \label{4.12}
    \int_\xi d\omega
    =\int_{\partial\xi}\omega+(-1)^{\ell}~d\int_\xi\omega
    \qquad\hbox{for $\omega\in\Omega^n(|N\U'|\times|N\U|)$}.
  \end{align}
  We refer to \cite{DL} for further details on fibre integration of
  simplicial forms.
\end{remark}


\medbreak
\section{Secondary characteristic classes}\label{five}

In this section we reformulate the classical constructions of
secondary characteristic classes and `characters' for connections on 
principal $G$-bundles in terms of simplicial forms. For the classical
constructions we refer to Kamber--Tondeur~\cite{KT}, Chern--Simons~\cite{CnS}, 
Cheeger--Simons~\cite{CrS} or Dupont--Kamber~\cite{DK}.

In the following $p\colon P\to X$ is a smooth principal $G$-bundle, $G$
a Lie-group with only finitely many components and $K\subseteq G$ is
the maximal compact subgroup. As in section~\ref{one} we fix an invariant
homogeneous polynomial $Q\in I^{n+1}(G)$, $n\geq0$, such that one of
the following 2 cases occur:

\paritem{Case I:} $Q\in\ker(I^{n+1}(G)\to I^{n+1}(K))$.

\paritem{Case II:} $Q\in I^{n+1}_\Z(G)$, that is, there exists an
integral class $u\in H^{2n+2}(BK,\Z)$ representing the Chern-Weil
image of $Q$ in $H^*(BG,\R)\cong H^*(BK,\R)$.

Let us introduce the notation
\begin{align}
  \label{5.1}
  H^{\ell+1}_\D(X)=\Omega^\ell (X)/ d\Omega^{\ell-1}(X)
\end{align}
In the notation of ~\cite{B2}, $H^{\ell+1}_\D(X)$ is canonically 
isomorphic to the smooth Deligne cohomology group 
$H^{\ell+1}_{\D,\infty} (X, 0(\ell+1))$, 
with '0' denoting the 0-ring, that is the hypercohomology group 
$\mathbf{H}^\ell (X, \underline{\Omega}^{(\ell)})$ 
with values in the truncated sheaf complex 
$$
\underline{\Omega}^{(\ell)} ~\colon~ 
\underline{\Omega}^0 \ra \underline{\Omega}^1 \ra \ldots 
\ra \underline{\Omega}^\ell~. 
$$ 
The elements $[\omega] \in H^{\ell+1}_\D(X)$ can be interpreted as equivalence 
classes of connections on the trivial $\ell$--gerbe $\theta = 0$ by setting 
\begin{equation}
  \label{5.1a}
  \omega^0=\veps^* \omega~, \qquad 
  F_{\omega}= d\omega~, \qquad 
  \check{\delta} \omega^0=0~, \qquad 
\omega^1=\ldots =\omega^{\ell} =0.
\end{equation} 
Clearly the connection is flat if and only if $F_\omega = d\omega=0$, 
that is $[\omega] \in H^{\ell} (X,\R)$. Com\-bi\-ning this with \eqref{2.11}, 
the data in \eqref{5.1a} determine a commutative diagram with exact rows 
\begin{equation}
   \label{5.1b}
   \xymatrix{
0 \ar[r] &H^\ell (X,\R) \ar[d] \ar[r]& H_\D^{\ell+1} (X) \ar[d] \ar[r] &
\Omega^\ell (X)/ \Omega_{\rm cl}^{\ell} (X) \ar[r] \ar[d]^d  &0 \\  
0 \ar[r] &H^\ell(X,\R/\Z) \ar[r] &H_\D^{\ell+1}(X,\Z) \ar[r]^{d_*} 
&\Omega_{\rm cl}^{\ell+1}(X,\Z) \ar[r] &0. }
\end{equation}
Further, it is easy to see that the center vertical arrow in diagram 
\eqref{5.1b} is induced by the exact hypercohomolgy sequence associated 
to the exact triangle of complexes
\begin{align*}
    \xymatrix{
      {\underline{\Omega}^{(\ell)}[-1]}
      \ar[rr]&&
      {\{\Z \ra \underline{\Omega}^0 \ra \underline{\Omega}^1 
      \ra \ldots \ra \underline{\Omega}^\ell\}} \ar[dl]\\
      &{\Z}
      \ar[ul]^{+1} }. 
\end{align*}

\medbreak
With this notation the \emph{secondary characteristic class}
associated to $Q$ (case I) or $(Q,u)$ (case II) for a connection $A$
on $P\to X$ is a class
\begin{equation}
\label{5.2}
  \begin{aligned}{}
    [\Lambda(Q,A)]&\in H^{2n+2}_\D(X)
    \qquad\hbox{in case I,}\\
    [\Lambda(Q,u,A)]&\in H^{2n+2}_\D(X,\Z)
    \qquad\hbox{in case II.}\\
  \end{aligned}
\end{equation}
Note that that the characteristic classes in $H^*_\D(X)$ are defined by global 
forms, whereas the classes in $H^*_\D(X,\Z)$ are defined by simplicial forms. 

\medbreak
For the construction we need the following well-known lemma which we
include for completeness:

\begin{lemma}
  \label{5.3}
  Given a $G$-bundle $p\colon P\to X$ with connection $A$ and an
  integer $N$, there is a bundle map 
  \begin{align}
    \label{5.4}
    \xymatrix{
      {P}\ar[r]^-{\bar\psi}\ar[d]&
      {\bar P}\ar[d]\\
      {X}\ar[r]^-\psi&
      {\bar X}}
  \end{align}
  and a connection $\bar A$ on $\bar P$ such that $\bar P$ is
  $N$-connected and such that $A=\bar\psi^*\bar A$.
\end{lemma}

\begin{proof}
  By choosing $\bar X$ to be a smooth approximation to the classifying
  space $BG$ we can clearly establish the bundle map in~\eqref{5.4}
  with $\bar P$ $N$-connected. Furthermore, by multiplying $\bar P$
  with a Euclidean space, the classifying map $\psi$ can be assumed to
  be an embedding. Then the connection $A$ on $P$ clearly extends over
  a tubular neighborhood of $X$ on $\bar X$ and subsequently over all
  of $\bar X$ by use of a partition of unity.
\end{proof}

\begin{remarks}
\label{rem4}
  \paritem{1.} Since the classifying map $\bar X\to BG$ for $\bar P$
  is unique up to homotopy, we have a natural identification of the
  cohomology $H^k(\bar X,\Z)\cong H^k(BG,\Z)$ for $k\leq N$.
  
  \paritem{2.} There is a functorial construction of the bundle map
  in~\eqref{5.4} using simplicial manifolds which however requires the
  use of multi-simplicial constructions for the Deligne cohomology
  (cf. \cite{DK}, \cite{DHZ}).
\end{remarks}

The classes in~\eqref{5.2} are now constructed as follows: Choose a
bundle map and connection $\bar A$ as in Lemma~\ref{5.3} with $N>2n+2$
and choose compatible good coverings $\U=\{U_i\}_{i\in I}$ and
$\bar\U=\{\bar U_{\bar\imath}\}_{\bar\imath\in\bar I}$ of $X$ respectively
$\bar X$. Also in case~II choose a representative
$\bar\gamma\in\Omega_\Z^{2n+2}(|N\U|)$ for the cohomology class 
$u\in H^{2n+2}(|N \bar{\U}|,\Z)\cong H^{2n+2}(BG,\Z)$. Then for $F_A$ and 
$F_{\bar A}$ the curvature forms for $A$ respectively $\bar A$, 
we can find (normal simplicial) forms $\Lambda(Q,\bar A)$ respectively
$\Lambda(Q,u,\bar A)$ such that
\begin{equation}
\label{5.5}
  \begin{aligned}
    Q(F_{\bar A}^{n+1})&=d\Lambda(Q,\bar A)
    \qquad\hbox{in case I},\\
    \epsilon^*Q(F_{\bar A}^{n+1})-\bar\gamma&=d\Lambda(Q,u,\bar A)
    \qquad\hbox{in case II},
  \end{aligned}
\end{equation}
and we put
\begin{equation}
\label{5.6}
  \begin{aligned}
    \Lambda(Q,A)&=\psi^*\Lambda(Q,\bar A)\in\Omega^{2n+1}(X)
    \qquad\hbox{in case I},\\
    \Lambda(Q,u,A)&=\psi^*\Lambda(Q,u,\bar A)\in\Omega_{\R/\Z}^{2n+1}(|N\U|)
    \qquad\hbox{in case II}.
  \end{aligned}
\end{equation}

\begin{proposition}
  \label{5.7}
  \paritem{$(i)$} The classes $[\Lambda(Q,A)]$, respectively
  $[\Lambda(Q,u,A)]$ in~\eqref{5.2} are well-defined given $\bar P$
  and $\bar A$.

  \paritem{$(ii)$} They are independent of the choice of $\bar P$ 
  and $\bar A$.

  \paritem{$(iii)$} They are natural with respect to bundle maps and
  compatible coverings.

  \paritem{$(iv)$} Curvature formula $:$~ 
  \begin{equation}\label{5.8}
    \begin{aligned}
      d\Lambda(Q,A)&=Q(F_A^{n+1})
      \qquad\hbox{in case I}\\
      d\Lambda(Q,u,A)&=\epsilon^*Q(F_A^{n+1})-\gamma
      \qquad\hbox{in case II}
    \end{aligned}
  \end{equation}
  where $\gamma=\psi^*\bar\gamma\in\Omega_\Z(|N\U|)$ represents the
  characteristic class $u(P)$ associated with $u$.

  \paritem{$(v)$} If $Q(F_A^{n+1})=0$, then 
  \begin{equation}\label{5.9}
    \begin{aligned}{}
      [~\Lambda(Q,A)~]&\in H^{2n+1}(X,\R)
      \qquad\hbox{in case I}\\
      [~\Lambda(Q,u,A)~]&\in H^{2n+1}(X,\R/\Z)
      \qquad\hbox{in case II},
    \end{aligned}
  \end{equation}  
and
  \begin{equation}
    \label{5.10}
    d_*[\Lambda(Q,u,A)]=-u(P)
  \end{equation}
  where $d_*\colon H^{2n+1}(X,\R/\Z)\to H^{2n+2}(X,\Z)$ is the
  Bockstein homomorphism.
\end{proposition}

\begin{proof}
  (i), (iii), (iv), and~(v) are obvious from the construction
  in~\eqref{5.5} and~\eqref{5.6}. Finally for (ii), let 
  $\bar\psi'\colon P\to\bar P'$ and $\bar A'$ be another choice 
  of bundle map and connection as in Lemma~\ref{5.3}. Then
  \begin{align*}
    \xymatrix{
      {P}\ar[d]\ar[r]^-{\bar\psi\times\bar\psi'}&
      {\bar P\times\bar P'}\ar[d]\\
      {X}\ar[r]&
      {(\bar P\times\bar P')/G}}
  \end{align*}
  is also a bundle map of the required form and 
  $A_t=(1-t)\bar A+t\bar A'~, ~t\in[0,1]$ gives a family of connections 
  on $\bar P\times\bar P'$ pulling back to the constant family $A$ in $P$. 
  The claim therefore follows from the following more general formula
  (with $\dd{A_t}t=0$).
\end{proof}

\begin{lemma}
  \label{5.11} Variational formula $:$ 
  Let $A_t$, $t\in[0,1]$, be a smooth family 
  of connections on $P\to X$ and let $\wti A$ denote the corresponding
  connection on $P\times[0,1]$ over $X\times[0,1]$. Then we have on
  $\Omega^{2n+1}(X)$ respectively $\Omega^{2n+1}(|N\U|):$
\begin{equation}
\label{5.12}
\begin{aligned}
      \Lambda(Q,A_1)-\Lambda(Q,A_0)& = 
(n+1)\int_0^1Q\Big(\dd{A_t}{t} \wedge F_{A_t}^n\Big)\,dt
      +d\int_0^1 i_{\dd{}{t}}\Lambda(Q,\wti A)\,dt~, 
      \\
\Lambda(Q,u,A_1)-\Lambda(Q,u,A_0)& = 
      (n+1) \epsilon^* \int_0^1 Q\Big(\dd{A_t}{t} \wedge F_{A_t}^n\Big) dt 
      +d\int_0^1 i_{\dd{}{t}}\Lambda(Q,u,\wti A) dt
\end{aligned}
\end{equation}
  in cases~I and~II respectively.
\end{lemma}

\begin{proof}
  Notice that the connection $\wti A$ on $P\times I$ satisfies
  $i_{\dd{}t}\wti A=0$. Hence for the curvature $F_{\wti A}=
  d\wti A+\frac{1}{2}[\wti A,\wti A]$ we have
  \begin{align*}
    i_{\dd{}t}F_{\wti A}=i_{\dd{}t}d\wti A=\dd{A_t}{t}.
  \end{align*}
  In case~II say, we therefore obtain from~\eqref{5.5}:
  \begin{equation}
  \label{5.13}
    \begin{aligned}
      \dd{}t\Lambda(Q,u,A_t)-di_{\dd{}{t}}\Lambda(Q,u,\wti A)
      &=i_{\dd{}{t}} d\Lambda(Q,u,\wti A)
      \\&=\epsilon^*i_{\dd{}{t}}Q(F_{\wti A}^{n+1})
      \\&=(n+1) ~\epsilon^* ~Q\Big(\dd{A_t}{t} \wedge F_{A_t}^n\Big),
    \end{aligned}
  \end{equation}  
since we can choose the representing integral form for $u$
  independent of $t$. Formula~\eqref{5.12} now follows
  from~\eqref{5.13} by integration.
\end{proof}

\medbreak
The invariants in \eqref{5.6} have certain 
\emph{multi\-pli\-cative pro\-per\-ties} which we state next. 

\begin{proposition}
  \label{5.18}
  \paritem{$(i)$} 
  For $Q_1$ and $Q_2$ both satisfying case I, we have
\begin{equation}
  \label{5.19}
    \begin{aligned}{}
    [~\Lambda(Q_1 Q_2, A)~] &= [~Q_1(A) \wedge \Lambda(Q_2,A)~] \\
    &= [~\Lambda(Q_1,A) \wedge Q_2(A)~] \in H_\D^* (X). 
    \end{aligned}
\end{equation}  

  \paritem{$(ii)$} 
  In case II, let $u_1, u_2$ and $u_1\cup u_2 \in H^*(BG,\Z)$ be 
  represented by integral forms $\gamma_1, \gamma_2$ and $\gamma_3$ 
  respectively, and choose the form $\mu$ such that 
  $d\mu = \gamma_1\wedge\gamma_2 - \gamma_3$.
  Then we have in $H_\D^* (X, \Z)~:$ 
\begin{equation*}
    \begin{aligned}{}
    [~\Lambda(Q_1Q_2,u_1\cup u_2,A)~] &= 
    [~\Lambda(Q_1,u_1,A) \wedge \psi^*\gamma_2 + 
    \epsilon^* Q_1(A) \wedge \Lambda(Q_2,u_2,A) - \psi^*\mu~] \\
    &= [~\psi^*\gamma_1 \wedge \Lambda(Q_2,u_2,A) + 
    \Lambda(Q_1,u_1,A)\wedge \epsilon^* Q_2 (A) - \psi^*\mu~].
    \end{aligned}
\end{equation*}  
\end{proposition}

\begin{proof}
This is straightforward from the definitions in \eqref{5.6}. 
\end{proof}

\medbreak
We now apply proposition~\ref{5.7} to the case of foliated 
bundles in the sense of Kamber--Tondeur~\cite{KT}. We recall that a
principal $G$-bundle $p\colon P\to X$ is \emph{foliated} if there 
are given two foliations $\wbar\F$ on $P$, $\F$ on $X$ such that
\newdimen\banjo     
\banjo=\linewidth   
\advance\banjo-4em  
  \begin{align}
  \label{5.14}
  \vcenter{
    \hbox{\parbox\banjo{\paritem{(i)} $\wbar\F$ is given by a 
    $G$--equivariant involutive subbundle $T \wbar\F \subset T P$, 
    that is the action by $G$ on $P$ permutes the leaves of $\bar\F$,}}
    \hbox{\parbox\banjo{\paritem{(ii)} for each $u\in P$ the
    differential $p_*\colon T_u\wbar\F\to T_{p(u)}\F$ is an
    isomorphism.}}}
\end{align}
Also the \emph{codimension} of the foliated bundle is by definition
the codimension of $\F$ in $X$. It is well-known that a foliated
$G$-bundle $p\colon P\to X$ has an \emph{adapted connection}, i.e., 
a connection $A$ satisfying $A(v)=0$ for $v\in T_u\wbar\F$, $u\in P$. 
Then it follows that the curvature form $F_A$ satisfies 
$F_A \in J$, where $J$ is the defining ideal of the foliation $\F$. 
For the codimension $q$ of $\F$, we have $J^{q+1} = 0$ and 
the curvature form satisfies $F_A^{q+1}\equiv 0$.

\begin{theorem}
  \label{5.15}
  \paritem{$(i)$} The classes $[\Lambda(Q,A)]$ respectively
  $[\Lambda(Q,u,A)]$ in~\eqref{5.2} are well-defined. 

  \paritem{$(ii)$} They are natural with respect to maps 
  of foliated bundles.

  \paritem{$(iii)$} Curvature formula $:$ ~We have 
  \begin{equation}
  \label{5.16}
    \begin{aligned}
      d\Lambda(Q,A)&=Q(F_A^{n+1})
      \qquad\hbox{in case I}\\
      d\Lambda(Q,u,A)&=\epsilon^*Q(F_A^{n+1})-\gamma
      \qquad\hbox{in case II}
    \end{aligned}
  \end{equation}
  where $\gamma=\psi^*\bar\gamma\in\Omega_\Z(|N\U|)$ represents the
  characteristic class $u(P)$ associated with $u \in H^{2n+2} (BK, \Z)$. 

  \paritem{$(iv)$} If $~n \geq q$, then $Q(F_A^{n+1})\in J^{q+1}=0$ and 
  \begin{equation}
  \label{5.17}
    \begin{aligned}{}
      [~\Lambda(Q,A)~]&\in H^{2n+1}(X,\R)
      \qquad\hbox{in case I}\\
      [~\Lambda(Q,u,A)~]&\in H^{2n+1}(X,\R/\Z)
      \qquad\hbox{in case II}. 
    \end{aligned}
  \end{equation}  
Moreover these classes are independent of the choice of adapted 
connection $A$. 

  \paritem{$(v)$} Rigidity $:$ ~
  If $n \geq q+1$, then the cohomology classes in 
  $(iv)$ are rigid under variation of the foliated structure 
  $(P, \wbar\F) \to (X,\F)$. 
\end{theorem}

\begin{proof}
  (i) to (iii) follow from the construction
  in~\eqref{5.5}, \eqref{5.6} and from proposition~\ref{5.7}. 
  The statements in (iv) and (v) essentially follow from the 
  variational formulas in \eqref{5.12}. 
  \eqref{5.17} in (iv) follows directly from \eqref{5.16}. 
  For the last statement in (iv), let $A'$ be another choice for 
  the adapted connection. 
  Then the family of adapted connections $A_t$ given by the convex 
  combination $A_t=(1-t)A+tA'~, ~t\in[0,1]$ satisfies 
  $\dd{A_t}{t}=A'-A=\alpha\in J$. 
  Thus we have $Q(\alpha\wedge F_{A_t}^{n})\in J^{q+1}=0$ for 
  $n \geq q$ and the statement follows from \eqref{5.12}. 
  For (v), let $(\wbar\F_t, \F_t)~, ~t\in[0,1]$, be a smooth family 
  of foliated structures on $P\to X$. 
  Let $A_t$, $t\in[0,1]$, be a smooth family 
  of $(\wbar\F_t, \F_t)$--adapted connections on $P\to X$. 
  Then for $n\geq q+1$ we have 
  $Q(\dd{A_t}{t}\wedge F_{A_t}^{n})\in J_t^{q+1}=0$ and (v) follows 
  also from \eqref{5.12}. 
\end{proof}

\begin{remarks}
\label{rem5}
  \paritem{1.} 
  Theorem~\ref{5.15} is essentially a reformulation 
  of theorem~2.2 in~\cite{DK}. 
  The above constructions could of course be extended to
  define more general characteristic classes associated to elements 
  in (the cohomology of) the relative Weil algebra 
  $F^{2(q+1)} W (G, K)$ as in~\cite{DK}.

  \paritem{2.} 
  Following Kamber-Tondeur \cite{KT},~section~2.24 we call the adapted 
  connection $A$ \emph{basic} if the Lie derivative $L_X A=i_X dA$ vanishes
  for all $\wbar\F$-horizontal vector fields $X$ on $P$ or equivalently
  if $i_X F_A=0$, that is $F_A \in J^2$. 
  If we can choose the connections in Theorem~\ref{5.15} to be basic, 
  then the condition $n\geq q$ in (iv) can be replaced by $2n\geq q$ 
  and the condition $n\geq q+1$ in (v) can be replaced by $2n \geq q+1$. 
  In fact, we have $Q(\alpha\wedge F_{A_t}^{n}) \in J^{2n+1}$ 
  in (iv), and $Q(\dd{A_t}{t}\wedge F_{A_t}^{n}) \in J_t^{2n}$ 
  in (v). 
\end{remarks} 


\medbreak
\section{Invariants for families of connections}\label{six}

We now return to the situation of a family of principal $G$-bundles
with connections as in Definition~\ref{1.1}. That is, (i) $\pi\colon
Y\to Z$ is a $\Diff^+(X)$-fibre bundle with fibre $X$, (ii) $p\colon
E\to Y$ is a principal $G$-bundle, and (iii) $A=\{A_z\mid z\in Z\}$ is
a family of connections on $P_z=E~| {X_z}~, ~z\in Z$, where 
$X_z = \pi\inv(z)$. 
Also $\V=\{V_j\}_{j\in J}$ and $\U=\{U_i\}_{i\in I}$ are good 
coverings of $Y$ respectively $Z$. Finally $Q\in I^{n+1}(G)$ is an 
invariant polynomial satisfying case~I or~II as in section~\ref{five}. 
Our main result is the following:

\begin{theorem}
  \label{6.1}
  Suppose $\partial X=\emptyset$ and $\dim X=2n+1-\ell$,
  $0\leq\ell\leq2n+1$. Also let $B$ be a global connection on $E$
  extending the family $A$. Then the following holds:
  
  \paritem{$(i)$} The $(simplicial)$ $\ell$-form defined by
  \begin{equation}\label{6.2}
    \begin{aligned}
      \Lambda_{Y/Z}(Q,B)&=\int_{Y/Z}\Lambda(Q,B)
      \qquad\hbox{in case I}\\
      \Lambda_{Y/Z}(Q,u,B)&=\int_{Y/Z}\Lambda(Q,u,B)
      \qquad\hbox{in case II}
    \end{aligned}
  \end{equation}  
  gives well-defined classes in $H_\D^{\ell+1}(Z)$, respectively
  $H_\D^{\ell+1}(Z,\Z)$, functorial with respect to bundle maps
  \begin{align}
    \label{6.3}
    \vcenter{\xymatrix{
        E'\ar[d]\ar[r]&E\ar[d]\\
        Y'\ar[d]\ar[r]&Y\ar[d]\\
        Z'\ar[r]      &Z}}
  \end{align}
  and the induced connections.
  
  \paritem{$(ii)$} These classes are independent of the 
  choice of the global extension $B$ provided that 
  $F_{A_z}^{n+1-\ell}=0$ for all $z\in Z$. 

  \paritem{$(iii)$} Curvature formula $:$ We have 
  \begin{equation}
  \label{6.4}
    \begin{aligned}
      \int_{Y/Z}Q(F_B^{n+1})&=(-1)^{\ell-1}~d\Lambda_{Y/Z}(Q,B)
      \qquad\hbox{in case I}\\
      \epsilon^*\int_{Y/Z}Q(F_B^{n+1})-\int_{Y/Z}\gamma
      &=(-1)^{\ell-1}~d\Lambda_{Y/Z}(Q,u,B)
      \qquad\hbox{in case II}
    \end{aligned}
  \end{equation}  
where $\gamma$ represents $u(E)\in H^{2n+2}(Y,\Z)$.

  \paritem{$(iv)$} In particular in case~II we have in
  $H^{\ell+1}(Z,\Z) :$
  \begin{align*}
    d_*~[~\Lambda_{Y/Z}(Q,u,B)~]=(-1)^{\ell}~\pi_!(u(E)).
  \end{align*}
  
  \paritem{$(v)$} If $F_{A_z}^{n-\ell}=0$ for all $z\in Z$ then
  $\Lambda_{Y/Z} (Q,B)$, respectively $\Lambda_{Y/Z} (Q,u,B)$ 
  are closed, respectively closed $\mod \Z$, and 
  \begin{equation}
  \label{6.5}
    \begin{aligned}{}
      [~\Lambda_{Y/Z}(Q,B)~]&\in H^\ell(Z,\R)
      \qquad\hbox{in case I}\\
      [~\Lambda_{Y/Z}(Q,u,B)~]&\in H^\ell(Z,\R/\Z)
      \qquad\hbox{in case II}
    \end{aligned}
  \end{equation}  
are well-defined invariants of the family $\{A_z\mid z\in Z\}$.
\end{theorem}

\begin{proof}
  Again (i), (iii), and (iv) follow from the definitions and the
  properties of fibre integration listed in proposition~\ref{4.10}.
  Also (ii) follows from Lemma~\ref{5.11} applied to the family 
  $B_t=(1-t)B_0+tB_1$, $t\in[0,1]$ for $B_0,B_1$ two choices of 
  global connections extending the family $A$.
  The first statements 
  in (v) are a consequence of formula~\eqref{6.4} and the 
  properties of fibre integration. 
  Finally, the last statement in (v) follows from (ii), since 
  the curvature assumption in (v) is stronger than the assumption 
  in (ii). 
\end{proof}

\begin{remarks}
\label{rem6}
  \paritem{1.} Theorems~\ref{1.2} and~\ref{1.3} are reformulations of
  Theorem~\ref{6.1}. In fact the $\ell$--gerbe $\theta=\theta(Q,u,B)$
  with connection $\omega=(\omega^0,\ldots,\omega^\ell)$ in
  Theorem~\ref{1.3} is given by the formulas in~\eqref{3.15} for
  $\Lambda=\Lambda_{Y/Z}(Q,u,B)$. 
  
  \paritem{2.} In particular we recover from Theorem~\ref{6.1} the
  construction of the Quillen `determinant line bundle' and their
  Hermitian connections as in \cite{DJ} by taking $\ell=1$ and 
  specializing $Y$ to a product and $\pi\colon X\times Z\to Z$ 
  the projection on $Z$ (compare also example \ref{7.20}).
  
  \paritem{3.} Again in the product situation $\pi\colon X\times Z\to
  Z$ and the covering $\V=\{U'_j\times U_i\}_{(i,j)\in I\times J}$ as
  in remark~\ref{rem3}, 
  we can define more generally for $\xi\in C^\Sing_{2n+1-\ell}(X)$ 
  respectively $\xi\in \check C^\Sing_{2n+1-\ell}(N\U')$ the invariant
  \begin{equation}
  \label{6.6}
    \begin{aligned}
      \Lambda_\xi(Q,B)&=\int_\xi\Lambda(Q,B)
      \qquad\hbox{in case I}\\
      \Lambda_\xi(Q,u,B)&=\int_\xi\Lambda(Q,u,B)
      \qquad\hbox{in case II}
    \end{aligned}
  \end{equation}  
Then we obtain from~\eqref{4.12} and~\eqref{5.8}:
  \begin{equation}
  \label{6.7}
    \begin{aligned}
      \int_\xi Q(F_B^{n+1})
      &=\Lambda_{\partial\xi}(Q,B) + (-1)^{\ell-1}~d\Lambda_\xi(Q,B)
      \qquad\hbox{in case I}\\
      \epsilon^*\int_\xi Q(F_B^{n+1}) - \int_\xi\gamma
      &=\Lambda_{\partial\xi}(Q,u,B) + (-1)^{\ell-1}~d\Lambda_\xi(Q,u,B)
      \qquad\hbox{in case II}.
    \end{aligned}
  \end{equation}
Here $\gamma$ represents $u(E)$ in $\Omega_\Z^{2n+2}(|N\U|)$, hence
  in particular $\int_\xi\gamma$ is integral. Thus, under the
  appropriate vanishing conditions for the fibre curvature the left
  hand side of~\eqref{6.7} is going to vanish (mod $\Z$ in case~II). 
  Hence $\Lambda_\xi(Q,B)$, respectively $\Lambda_\xi(Q,u,B)$, defines 
  a cycle in the total complex of the bicomplex 
\begin{equation*}
    \begin{aligned}
      \Hom (C_* (X) &, \Omega^* (Z))  
      \qquad\hbox{in case I}\\
      \Hom (\check{C}_* (\U') &, \Omega^*_{\R/\Z}(||N\U||))
      \qquad\hbox{in case II}.
    \end{aligned}
\end{equation*}  
  Notice that for $\ell=0$, 
  $\Lambda_\xi(Q,u,B)$ is essentially 
  the `Chern--Simons section' of the line bundle given by
  $\Lambda_{\partial\xi}(Q,u,B)$ as defined in \cite{DJ}.
\end{remarks}

\medbreak
We can now apply theorem~\ref{6.1} to the general case of families 
of foliated bundles. By a \emph{family} of \emph{foliated} $G$-bundles 
of codimension $q$ we mean the following:
\begin{align}
  \label{6.8}
  \vcenter{
    \hbox{\parbox\banjo{\paritem{(i)} $\pi\colon Y\to Z$ is a
    $\Diff^+(X)$-fibre bundle with fibre $X$.}}
    \hbox{\parbox\banjo{\paritem{(ii)} $p\colon E\to Y$ is a principal
        $G$-bundle.}}
    \hbox{\parbox\banjo{\paritem{(iii)} $\wbar\F, ~\F$ are foliations 
    of $E$, respectively $Y$, such that $T \F \subset T (\pi)$, 
    respectively $T \wbar\F \subset T (\pi\circ p)$ are involutive 
    ($G$--equivariant) subbundles, inducing foliated structures 
    $(\wbar\F_z, \F_z)$ of codimension $q$ in the principal bundles 
    $p_z : P_z \to X_z$ 
    for $z\in Z$.}}}
\end{align}
In this situation $(\wbar\F, \F)$ makes $p\colon E\to Y$ into a foliated 
$G$-bundle. By restriction to $T (\pi\circ p) \subset T E$, a global 
adapted connection $B$ induces a smooth family $A = \{ A_z \}$ of adapted 
connections on the principal bundles $p_z : P_z \to X_z~, ~z \in Z$, 
satisfying the curvature condition $F_{A_z}^{q+1}  = 0$. 
Conversely, any global extension $B$ of a smooth family $A = \{ A_z \}$ 
of adapted connections is adapted to $(\wbar\F,\F)$. 
Thus by choosing a global adapted connection $B$, we conclude the 
following from Theorem \ref{6.1}: 

\begin{theorem}
  \label{6.9}
  Suppose $\partial X=\emptyset$ and $\dim X=2n+1-\ell$,
  $0\leq\ell\leq2n+1$. Let $B$ be an adapted connection for the 
  family of foliated bundles of codimension $q$ as above. 
  Then the following holds ~$:$ 
  
  \paritem{$(i)$} The classes 
  \begin{equation}
  \label{6.9a}
    \begin{aligned}{}
      [~\Lambda_{Y/Z}(Q,B)~] &\in H_\D^{\ell+1}(Z)
      \qquad\hbox{in case I}\\
      [~\Lambda_{Y/Z}(Q,u,B)~] &\in H_\D^{\ell+1}(Z,\Z)
      \qquad\hbox{in case II}
    \end{aligned}
  \end{equation}  
  are well-defined and independent of the choice of adapted
  connection $B$ if $n-\ell\geq q$. 

  \paritem{$(ii)$} Suppose that $n-\ell>q$. Then $\Lambda_{Y/Z} (Q,B)$, 
  respectively $\Lambda_{Y/Z} (Q,u,B)$ are closed, re\-spect\-ively 
  closed $\mod \Z$ and
  \begin{equation}
  \label{6.10}
    \begin{aligned}{}
      [~\Lambda_{Y/Z}(Q,B)~]&\in H^\ell(Z,\R)
      \qquad\hbox{in case I}\\
      [~\Lambda_{Y/Z}(Q,u,B)~]&\in H^\ell(Z,\R/\Z)
      \qquad\hbox{in case II}
    \end{aligned}
  \end{equation} 
  are well-defined invariants of the family of foliated bundles. 

  \paritem{$(iii)$} Suppose again that $n-\ell>q$. Then the cohomology 
  classes $[~\Lambda_{Y/Z} (Q,B)~]$, respectively $[~\Lambda_{Y/Z} (Q,u,B)~]$ 
  in $\eqref{6.10}$ above are rigid, that is they are invariant under 
  $($germs of$)$ smooth deformations of the data in $\ref{6.8}, (iii)$. 
\end{theorem}

\medbreak
  In either case, we call the invariants in \eqref{6.9a}, respectively 
  in \eqref{6.10} the \emph{characteristc $\ell$--gerbe}, respectively 
  the \emph{characteristc flat $\ell$--gerbe} of the family of foliated 
  bundles, associated to the pair $(Q, u)$. 

\begin{proof}
  (i) needs some elaboration, since the family $A$ of adapted 
  connections on $T (\pi\circ p)$ is now not fixed. 
  We want to show that (i) follows from the variational formulas 
  in \eqref{5.12}. 
  Let $A, A'$ be two families of adapted connections along the fibres 
  and consider corresponding global extensions $B, B'$ of $A, A'$. 
  Then the convex combination $B_t=(1-t)B+tB'~, ~t\in[0,1]$ is an 
  extension of the adapted connection $A_t=(1-t)A+tA'~, ~t\in [0,1]$ 
  on the fibres. Further $B_t$ satisfies 
  $\dd{B_t}{t}=B'-B=\beta$, where the $\wbar\F$--transversal $1$--form 
  $\beta$ on $Y$ is of the form $\beta = \alpha + \gamma$, with 
  $\alpha = \alpha^{1,0} = A'-A$ on $T (\pi)$ being fibrewise in the 
  ideal $J_z$ of $\F_z$, that is $\alpha$ vanishes on the subbundle 
  $T \wbar\F \subset T (\pi)$, and $\gamma = \gamma^{0,1}$ being of 
  type $(0,1)$ on $Y$, that is $\gamma$ vanishes on the subbundle 
  $T (\pi) \subset T Y$. 
  Thus we have, observing that $F_{B_t}^{2,0} = F_{A_t}$, 
\begin{equation}
  \label{6.11}
  \begin{aligned}
  \int_{\Delta^1} Q (F_{\wti{B}}^{n+1}) &= 
  (n+1) \int_{\Delta^1} dt \wedge Q(\tp{\beta}{F_{B_t}^{n}}) \\ 
  &= (n+1)\int_0^1 
  Q\Big( (\alpha+\gamma) \wedge 
  \left( F_{B_t}^{2,0}+F_{B_t}^{1,1}+F_{B_t}^{0,2} \right)^n \Big)\,dt \\
  &= (n+1)\int_0^1 
  Q \Big(\alpha^{1.0} \wedge 
  \left( F_{A_t}+F_{B_t}^{1,1}+F_{B_t}^{0,2} \right)^n \Big)\,dt \\
  &+ (n+1)\int_0^1 
  Q \Big( \gamma^{0,1} \wedge 
  \left( F_{A_t}+F_{B_t}^{1,1}+F_{B_t}^{0,2} \right)^n\Big)\,dt~. 
  \end{aligned}
  \end{equation}  
  As we will have to integrate over the fibre, only the components of 
  type $(2n+1-\ell, \ell)$ of the $(2n+1)$--form in the integrand can 
  contribute non--trivial terms. 
  Therefore the relevant terms in the first summand of \eqref{6.11} 
  must contain $\alpha \wedge F_{A_t}^k$ for $k \geq n-\ell\geq q$, 
  that is $k+1 \geq (n-\ell)+1 \geq q+1$,
  while the relevant terms in the second summand of \eqref{6.11} 
  must contain $F_{A_t}^k$ for $k \geq n-(\ell-1) = (n-\ell)+1\geq q+1$. 
  Since $\alpha \wedge F_{A_t}^q=0$ and $F_{A_t}^{q+1}=0$, it follows 
  that all the relevant terms vanish in either case. 
  Thus \eqref{6.11} vanishes under integration over the fibre 
  and (i) follows from \eqref{5.12}. 
  
  (iii) is proved by the same argument, with the following modifications. 
  Let $A_\sigma~, ~\sigma \in [0, 1]$ be families of connections 
  along the fibres, adapted to the data $(\wbar\F_\sigma, ~\F_\sigma)$ 
  and let $B_\sigma$ be corresponding global 
  extensions of $A_\sigma$. Then we can write again 
  $\dd{B_\sigma}{\sigma} = \alpha_\sigma^{1,0} + \gamma_\sigma^{0,1}$ 
  as above, except that the horizontal forms $\alpha_\sigma^{1,0}$ on $T (\pi)$ 
  do not necessarily satisfy the fibrewise condition of being in the 
  ideal $(J_\sigma)_z$ of $(\F_\sigma)_z$. However, the argument in the 
  proof of (ii) remains valid, since the relevant terms 
  $\alpha_\sigma \wedge F_{A_\sigma}^k$ satisfy now 
  $k \geq n-\ell\geq q+1$ and therefore vanish. 

  Of course, (ii) follows from theoren \ref{6.1} (v). Explicitly, 
  we have to show that the curvature term in \eqref{6.4} vanishes under 
  the assumption $n-\ell>q$. 
  Writing $F_{B} = F_{A}+F_{B}^{1,1}+F_{B}^{0,2}$ as above and expanding 
  $$
  Q (F_B^{n+1}) =  Q \left( (F_{A}+F_{B}^{1,1}+F_{B}^{0,2})^{n+1} \right), 
  $$
  the claim follows by a counting argument similar to the one above. 
\end{proof}

\begin{remarks}
\label{rem7}
  \paritem{1.} 
  In applying the variation formula \eqref{5.12} in lemma \ref{5.11} 
  in the proof of theorem \ref{6.9}~(i), we observe that the 
  adapted connection $A = B \mid T (\pi\circ p)$ is not fixed during a 
  variation $B_t$ of $B$, but we still have 
  $F_{A_t}^{q+1} = 0~, ~t \in [0,1].$ 

  \paritem{2.} 
  The results of theorem \ref{6.9} apply in particular to families 
  $p\colon E \to Y$ of \emph{flat} bundles, that is $A=\{ A_z \}$ 
  is a family of flat connections on the $G$--principal bundles 
  $P_z = E~| {X_z} \to X_z$ for $z \in Z$. In this case we have 
  $q=0$ and $T \F = T (\pi)$ and the relevant conditions are 
  $n \geq \ell$ in (i) and $n > \ell$ in (ii) and (iii). This case occurs 
  in all examples in section \ref{seven} except for the last example 
  \ref{7.30}. Families of flat bundles are also considered in ~\cite{GJK}. 

\paritem{3.} 
  The characteristic classes of foliated bundles in theorem \ref{5.15} 
  can be recovered from theorem \ref{6.9} by taking 
  $Z = \{\pt\}$ and $\ell=0$. 
  In this case, we have $B=A$ and the non--integrated classes 
  $[~\Lambda (Q,A)~] \in H^{2n+1} (X,\R)$, respectively 
  $[~\Lambda(Q,u,A)~] \in H^{2n+1} (X,\R/\Z)$ are well-defined 
  under the assumption $n \geq q$ as in theorem \ref{5.15}, (iv). 
  Moreover, the restriction $\dim X=2n+1$ is obviously not necessary. 
\end{remarks} 

\medbreak
Recall that the family of adapted connections $A$ is basic if 
the Lie derivative $L_X A=i_X dA$ vanishes for all $\wbar\F$-horizontal 
vector fields $X$ on $E$ or equivalently if $i_X F_A=0$. 

\begin{proposition}
  \label{6.12}
  If we can choose $A$ basic in Theorem~$\ref{6.9}$, then the conditions 
  $n-\ell\geq q$ in~$(i)$, respectively $n-\ell>q$ in~$(ii)$, 
  can be replaced by $2(n-\ell)\geq q$, respectively $2(n-\ell)>q$. 
\end{proposition}

\begin{proof}
  Again, we restrict attention the first statement (i) in theorem 
  \ref{6.9}. Counting powers of $J_t$ instead of curvature terms, 
  we see that the above estimates for the relevant terms in the proof 
  of (i) give $2k+1 \geq 2(n-\ell)+1$ 
  for the first summand of \eqref{6.11}, and 
  $2k \geq 2(n-(\ell-1))$, that is $2k \geq 2(n-\ell)+2$ 
  for the second summand of \eqref{6.11}. 
  Thus in either case, the condition $2(n-\ell)\geq q$ implies that the 
  relevant curvature terms are in $J^{q+1} = 0$.    
\end{proof}

\begin{remarks}
\label{rem8}
  \paritem{1.} 
  Observe that the global extensions $B_t$ on $E$, respectively 
  the connection $\wti{B}$ on $E\x[0,1]$ in the proof of 
  theorem \ref{6.9}~(i) will in general not be basic for the respective 
  foliated structures, even if $A, A'$ and hence $A_t$ are.  

  \paritem{2.} 
  One might expect the correct conditions in proposition \ref{6.12} 
  to be $n-\ell \geq q'$ in~(i), respectively $n-\ell > q'$ in~(i), 
  where $q' = [\frac{q}{2}]$, that is $q=2q'$ for $q$ even and 
  $q=2q'+1$ for $q$ odd. Then the basic vanishing property is 
  $F_A^{q'+1}=0$, since $2(q'+1) = 2q'+2 \geq q+1$. 
  However, for $q$ odd and $n-\ell = q'$, 
  the estimate $2k+1 \geq 2(n-\ell)+1$ gives 
  $2k+1 \geq 2(n-\ell)+1 = 2q'+1 = q$ which is not sufficient. 
\end{remarks}


\medbreak
\section{Examples}\label{seven}

In this section we give a few examples of increasing complexity.  

\begin{example}
  \label{7.1}
  We start with a simple example, found together with R.~Ljungmann,  
  which gives non--trivial classes in $H_\D^{2}(Z)$. 
  Let $X =T^2=\R^2/\Z^2$ and $Z=\R^2$ and 
  consider the trivial $\GL(1,\R)_+=\R_+^\x$--bundle $E$ over 
  $Y=T^2\times \R^2$ with coordinates $(x_1,x_2;z_1,z_2;\lambda)$. 
  If $\omega_0$ denotes the Cartan--Maurer form, then 
  $\omega = \omega_0 + B~, ~B = z_1 ~dx_1 + z_2 ~dx_2,$ defines a 
  foliated structure on $E$ which is flat along the fibres $T^2$ of 
  $\pi:Y \to Z$. In fact, the curvature on $Y$ is given by  
  $F = dB = dz_1\wedge dx_1 + dz_2\wedge dx_2$, 
  which is clearly of type $(1,1)$ and vanishes on every fibre 
  $T^2_z=\pi^{-1} (z)~, ~z=(z_1,z_2)$.  
  The flat structure of $E\mid T^2_z \to T^2_z$ is not trivial; 
  in fact, the holonomy depends on $z \in Z$ and is given by the homomorphism 
  $h_z : \Lambda \to \R_+^\x~, ~\Lambda = \pi_1(T^2_z) \cong \Z^2~,$ where 
\begin{align}
  \label{7.2}
  h_z (\lambda_1,\lambda_2) = e^{\lambda_1 z_1+\lambda_2 z_2}.
  \end{align}
  Since the Lie algebra of $\R_+^\times$ is $\R$ we can take the
  polynomial $Q(\xi)=\xi^2$ to obtain the $3$--form 
  \begin{align*}
  \Lambda(Q,B)=B\wedge dB= (dx_1\wedge dx_2) \wedge (z_2 ~dz_1-z_1 ~dz_2),
  \end{align*}
  with $d\Lambda(Q,B)=dB^2= -2 ~(dx_1\wedge dx_2) \wedge (dz_1\wedge dz_2)$. 
  Thus on $Z$ we have the characteristic form
  \begin{align}
  \label{7.3}
    \Lambda_{Y/Z}(Q,B)= \int_{T^2} B\wedge dB= z_2 ~dz_1-z_1 ~dz_2.
  \end{align}
  which defines a non-zero class in $H_\D^{2}(Z)$, and can be interpreted 
  as a connection in the trivial line bundle on $Z$ with curvature 
  $-2 ~dz_1\wedge dz_2 = -2 V$, where $V$ is the volume form on $Z=\R^2$. 
  
  Restricting $Y$ and \eqref{7.3} to $\S^1 \subset Z=\R^2$ by setting 
  $z_1=\cos\theta,~z_2=\sin\theta$, we obtain on $T^2\times \S^1$
  \begin{align*}
    \Lambda(Q,B)=B\wedge dB=-dx_1\wedge dx_2\wedge d\theta. 
  \end{align*} 
  Thus on $\S^1$ we have the characteristic form
  \begin{align}
  \label{7.4}
    \Lambda_{Y/\S^1}(Q,B)=-\Big(\int_{T^2} dx_1\wedge dx_2\Big)\,d\theta
    =-d\theta~, 
  \end{align}
  representing a non-zero element in $H^1(\S^1,\R)\cong \Hom(\Z, \R)=\R$.
  Thus the restriction of the class in \eqref{7.3} is closed, that is the 
  above line bundle is flat on $\S^1$ with holonomy determined by \eqref{7.4}. 
\end{example}

\begin{example}
  \label{7.5}
  More generally let $X=X_g$ be a surface of genus $g\geq2$ and let 
  $\{ \alpha_1, \beta_1,\ldots, \alpha_g, \beta_g \}$ be a set of 
  closed $1$--forms representing a symplectic basis for the 
  cup--product pairing in cohomology, that is 
\begin{equation*} 
\int_{X_g} \alpha_i\wedge \alpha_j=0 ~, \qquad 
\int_{X_g} \alpha_i\wedge \beta_j=\delta_{ij} ~, \qquad   
\int_{X_g} \beta_i\wedge \beta_j=0. 
\end{equation*}
  We let $Z=\R^{2g}$ with coordinates $(z_1,\ldots,z_{2g})$
  and again consider the foliated $\R_+^\x$--bundle $E$ with the 
  foliated structure given by the $1$--form $\omega = \omega_0 + B~,~
  B = z_1~\alpha_1+z_2~\beta_1 +\ldots+z_{2g-1}~\alpha_g+z_{2g}~\beta_g$, 
  similar to example \ref{7.1}. The curvature $F = dB$ on $Y$ is again   
  of type $(1,1)$ and vanishes on every fibre 
  $T^2_z=\pi^{-1} (z)~, ~z\in Z$.  
  The holonomy of the flat bundles $E~\vert X_{g, z}$ is determined 
  as a homomorphism $h_z : \Gamma \to \R_+^\x~, ~
  \Gamma = H_1 (X_g, \Z) \cong \Z^{2g}~,$ by a formula similar to \eqref{7.2}, 
  namely 
  \begin{align}
  \label{7.5a}
  h_z (\gamma_1,\ldots, \gamma_{2g}) = e^{\int_{\gamma_1}\alpha_1 +
  \int_{\gamma_2}\beta_1+ \ldots +\int_{\gamma_{2g-1}}\alpha_g +
  \int_{\gamma_{2g}}\beta_g}. 
  \end{align}
  Again we take the polynomial $Q(\xi)=\xi^2$ 
  to obtain the characteristic form $1$--form on $Z$: 
\begin{align}
  \label{7.6}
    \Lambda_{Y/Z}(Q,B)= \int_{X_g} B\wedge dB= 
    (z_2 ~dz_1-z_1 ~dz_2) + \ldots + (z_{2g} ~dz_{2g-1}-z_{2g-1} ~dz_{2g}), 
  \end{align}
  which defines a non-zero class in $H_\D^{2}(Z)$, and can be interpreted 
  as a connection in the trivial line bundle on $Z$ with curvature 
\begin{align}
  \label{7.7}
    d \Lambda_{Y/Z}(Q,B)= \int_{X_g} dB^2 = -2 \left(
    dz_1\wedge dz_2 + \ldots + dz_{2g-1}\wedge dz_{2g} \right), 
\end{align}
  Note that in this and the previous example we have $n=\ell=1$ and $q=0$. 
\end{example}

\begin{example}
  \label{7.8}
  This example is like example \ref{7.1}, but here we take 
  $X = T^k=\R^k/\Z^k$ and $Z=\R^k$ and consider again the trivial 
  $\GL(1,\R)_+=\R_+^\x$--bundle $E$ over $Y=T^k\times \R^k$ with 
  coordinates $(x_1,\ldots,x_k;z_1,\ldots,z_k;\lambda)$, with 
  the foliated structure given by the $1$--form $\omega = \omega_0 + B~, ~
  B=z_1 ~dx_1 + \ldots + z_k ~dx_k$. This foliated structure is flat along 
  the fibres $T^k$ of $\pi:Y \to Z$. In fact, we have for the curvature 
  $F = dB = dz_1\wedge dx_1 +\ldots + dz_k\wedge dx_k$, 
  which is of type $(1,1)$ and vanishes on every fibre 
  $T^k_z=\pi^{-1} (z)~, ~z=(z_1,\ldots,z_k)$.  
  As in example \ref{7.1}, the holonomy of the flat bundle 
  $E\mid T^k_z \to T^k_z$ depends on $z \in Z$ and is given by the 
  homomorphism 
  $h_z : \Lambda \to \R_+^\x~, ~\Lambda = \pi_1(T^k_z) \cong \Z^k~,$ where 
\begin{align}
  \label{7.9}
  h_z (\lambda_1,\ldots,\lambda_k) = e^{\lambda_1 z_1+\ldots+\lambda_k z_k}~, ~
  (\lambda_1,\ldots,\lambda_k) \in \Lambda~. 
\end{align}
  Now we take the polynomial $Q(\xi)=\xi^{n+1}~, ~k=n+1$ to obtain the 
  characteristic $(2n+1)$--form 
  \begin{align*}
  \Lambda(Q,B)&=B\wedge dB^n \\
  &= (-1)^{\binom{k}{2}}~(k-1)!~(dx_1\wedge\ldots \wedge dx_k) \wedge \\ 
  &~ \wedge \left( \sum_{j=1}^k ~(-1)^{j-1} ~
  z_j ~dz_1 \wedge \ldots \wedge \what{dz_j} 
  \ldots \wedge dz_k \right),
  \end{align*}
  Thus on $Z = \R^k$ we have the characteristic form
  \begin{equation}
  \begin{aligned}
  \label{7.10}
    \Lambda_{Y/Z}(Q,B)&= \int_{T^k} ~B \wedge dB^{k-1} \\ 
    &= (-1)^{\binom{k}{2}}~(k-1)!
    ~\sum_{j=1}^k ~(-1)^{j-1} ~z_j ~dz_1 \wedge \ldots \wedge \what{dz_j} 
    \ldots \wedge dz_k~, 
  \end{aligned}
  \end{equation}
  with curvature 
  \begin{align}
  \label{7.11}
  d\Lambda_{Y/Z}(Q,B) = (-1)^{\binom{k}{2}}~k! ~
  dz_1 \wedge \ldots \wedge dz_k = 
  (-1)^{\binom{k}{2}}~k! ~V~, 
  \end{align}
  with $V$ the volume form on $Z=\R^k$. 
  Hence \eqref{7.10} defines a non-zero class in $H_\D^{k}(Z)$. 
  
  Restricting $Y$ and \eqref{7.10} to 
  $\S^{k-1} = \{(z_1,\ldots,z_k)\mid \sum_{i=1}^k ~z_i^2=1\} 
  \subset Z=\R^k$, it is easy to see that $\Lambda_{Y/Z}(Q,B)$ 
  is a non--zero multiple of the volume form on $\S^{k-1}$ and 
  is clearly closed. 
  Thus we have $\Lambda_{Y/\S^{k-1}}(Q,B)\neq 0 \in H^{k-1}(\S^{k-1})$. 
  Note that in this example we have $k=n+1~, ~n=\ell$ and $q=0$. 
  We can interprete the invariants $\Lambda_{Y/Z}(Q,B)$, respectively 
  $\Lambda_{Y/\S^{k-1}}(Q,B)$ as (flat) connections on the trivial 
  $n=(k-1)$--gerbe as in \eqref{5.1a}. 

\end{example}

\medbreak 
So far, the examples have been for case I. 
The next two examples will be for case II.  

\begin{example}
  \label{7.15} The {\em Poincar\'e $(k-1)$--gerbe} 
  (cf. ~\cite{BMP},~\cite{GJK})~:~
This example is the case II analogue of example \ref{7.8}. 
Let $T$ be the $k$--dimensional real torus, that is
$T=\R^k/\Lambda$ for the rank $k$ integral lattice
$\Lambda\subset\R^k$~. The associated dual torus is defined as
\begin{equation}\label{7.16}
\what{T} = H^1(T, \R)/{H^1(T, \Z)} \ocong{\exp} 
{\Hom_{\Z}(\Lambda, \Un (1))} \cong \Un (1)^k~, 
\end{equation}
that is the points in $\what T$ parametrize flat unitary 
connections on the trivial line bundle
$\underline{\C} = T \times \C \lra T$, 
For $\xi \in \what{T}~,~ x \in \R^k~,~ a \in \Lambda$ and
$\lambda \in \C~,$ consider the equivalence relation
\begin{equation}
\label{7.17}
\begin{aligned}
\R^k \times \what{T} \times \C &\lra  \R^k \times \what{T}
\times \C/{\sim}~, \\ 
(x + a, \xi, \lambda) &\sim (x, \xi, \exp (2 \pi \iota~\xi(a)) \lambda)~. 
\end{aligned}
\end{equation}
The quotient space under `$\sim$' defines the Poincar\'e line
bundle $\P \lra T \x \what{T}$~. Let $\hat p$ denote the
projections of $T\times \what{T} \to \what{T}$. 
From \eqref{7.17} we see that the restriction 
$\P \mid \hat{p}^{-1}(\xi) \cong \L_{\xi}~,$ where the latter 
denotes the flat line bundle parametrized by $\xi \in \what{T}$. 
There exists a canonical unitary connection $B$ on the 
$\Un(1)$--principal bundle $p : E \to T\times \what{T}$ 
associated to $\P$, with curvature $F_B$ given by
\begin{equation} 
\label{7.18}
F_B = 2 \pi \iota~ \sum_{j=1}^k  d\xi^j \wedge dx_j~, 
\end{equation}
where $\{x_j\}$ are (flat) coordinates on $T$ and $\{\xi^j\}$ are
dual (flat) coordinates on $\what T$. $F_B$ is of type $(1,1)$ 
and therefore induces a family $A = \{ A_\xi \}$ of flat connections 
on the fibres $T \to T\times \what{T} \oset{\hat p} \what{T}$, 
Now we take $Q = C_1^{k}$ and $u=c_1^k$, where 
$c_1 \in H^2 (B\Un (1)),\Z)=H^2 (\C \PBbb^{\infty},\Z)\cong \Z [c_1]$ 
is the generator. 
Thus we obtain from theorem \ref{6.9} a (simplicial) characteristic 
$n=(k-1)$--gerbe 
\begin{equation}
\label{7.19}
[~\Lambda_{T \x \what{T}/\what{T}} (C_1^k, c_1^k, B)~] = \int_T ~
[~\Lambda (C_1^k, c_1^k, B)~] \in H_\D^{k} (\what{T}, \Z). 
\end{equation}

We remark that $T \x \what{T}$ has a canonical K\"ahler structure 
for which the Poincar\'e bundle $\P$ becomes a holomorphic 
line bundle such that $C_1 (\P) = \frac{1}{2 \pi \iota} F_B = \omega$. 
Here $\omega$ is the K\"ahler form, so that $T \x \what{T}$ has a Hodge 
structure. It follows that the curvature of the characteristic gerbe 
in \eqref{7.19} is a non--zero multiple of 
$\int_T ~\omega^k = V$, the volume form on $\what{T}$.
Note that in this and the previous example we have $n=\ell=k-1$ and $q=0$. 
\end{example}

\begin{example}
  \label{7.20} The {\em Quillen $1$--gerbe}~
  \cite{Q}, ~\cite{RSW}, ~\cite{DJ}~:~ 
This well-known complex line bundle with unitary connection associated 
to families of flat $\SU (2)$--bundles appears in our setup as a characteristic 
$1$--gerbe. We briefly recall this non--abelian example, referring to Ramadas--
Singer--Weitsman ~\cite{RSW}~ for details. 
Let $X=X_g$ be an oriented surface of genus $g~, ~G=\SU (2)$ and let $Z$  
be the smooth part of the representation variety $\Hom (\pi_1 (X_g), G) / G~.$
This is a symplectic manifold of dimension $6 (g-1)$ and the symplectic form 
is in fact the curvature form for the characteristic $1$--gerbe constructed below. 
The family $E \ra X_g \x Z$ is the tautological family of flat $\SU (2)$--bundles 
$P_\rho \to X_g$ determined by 
$\rho : \pi_1 (X_g) \to \SU (2)~, ~\rho \in Z~.$ 
The pair $(Q, u)$ is taken to be $Q = C_2~,$ the second Chern polynomial, and 
$u=c_2 \in H^4 (B\SU (2), \Z) \cong \Z [c_2]$ is the universal Chern class. 
Hence choosing a global $\SU (2)$--connection $B$ on $E$, extending the family $A$ 
of flat connections along the fibres $P_\rho \to X_g~, ~\rho \in Z$, we obtain 
from theorem \ref{6.9} the (simplicial) characteristic $1$--gerbe 
\begin{equation}
\label{7.21}
[~\Lambda_{X_g \x Z / Z} (C_2, c_2, B)~] = \int_{X_g} ~
[~\Lambda (C_2, c_2, B)~] \in H_\D^2 (Z, \Z). 
\end{equation}

\end{example}


\medbreak
The above examples are all cases where $q=0$, that is we have 
$T \F = T (\pi)$ and $A=\{ A_z \}$ is a family of 
\emph{flat} connections on the fibres $P_z \to X_z~, ~z\in Z$. 
We end with a case I example which relates to variations of the 
Godbillon--Vey invariant ~\cite{He1} and also gives some new classes 
of Godbillon--Vey type. 

\begin{example}
  \label{7.30} {\em Godbillon--Vey gerbes} for families of foliations~: ~
Let $\F$ be a family of transversally oriented foliations of codimension 
$q$ on $\pi:Y \to Z$ as in \eqref{6.8}, that is $T \F \subset T (\pi)$. 
The relative transversal bundle $Q_\F = T (\pi) / T \F$ has a natural 
foliated structure given by the partial Bott connection. 
On the oriented frame bundle $E = F_{\GL (q)^+} (Q_\F) \to Y$ this 
determines a foliated structure $\wbar\F$. We choose a family 
$A = \{ A_z \}$ of torsion--free, hence adapted connections along 
the fibres and extend it to a global connection $B$ on $E \to Y$. 
For given $n \geq q$, we consider invariant polynomials of the form 
$C_1 Q \in \ker (I (\GL (q, \R)^+) \to I (\SO (q))$, 
where $Q \in I^n (\GL (q, \R)^+)$ and 
$I (\GL (q, \R)^+) \cong \R [C_1, \ldots, C_q]$ 
is generated by the Chern polynomials $C_j$, that is the coefficients 
of $t^j$ in $\det (\Id+\frac{t}{2\pi} A)~, ~A \in \gl (q, \R)$.  
We have $C_1=\frac{1}{2\pi} ~\Tr$ and the kernel of the restriction to 
$I (\SO (q))$ is generated by the odd Chern polynomials $C_{2k+1}$.
Then $\Lambda (C_1 Q, B)$ is given by the $(2n+1)$--form 
\begin{equation}
\label{7.31}
\Lambda (C_1 Q, B)= \beta \wedge Q (F_B^n) 
\end{equation}
on $Y$, satisfying 
\begin{equation}
\label{7.32}
d \Lambda (C_1 Q, B) = 
d \beta \wedge Q (F_B^n) =
\frac{1}{2\pi} \Tr (F_B) \wedge Q (F_B^n) = 
C_1 (F_B) \wedge Q (F_B^n).  
\end{equation}
Here $\beta = \frac{1}{2\pi}s^* \Tr (B)$ is the pull--back of the 
trace of the connection form $B$ on $\det (E) = \Lambda^q (Q_\F)_0$ 
by a trivializing section $s : Y \to \Lambda^q (Q_\F)_0$ given by the 
transverse orientation on the normal bundle $Q_\F$. 
Note that $\beta$ satisfies $d \beta=C_1 (F_B)$ and that the choice 
$Q = C_1^n$ corresponds to the Godbillon--Vey form proper, that is 
$\Lambda (C_1^{n+1}, B)= \beta \wedge d\beta^n$. 
For $\ell$ satisfying $n-\ell \geq q$, \eqref{7.31} now gives rise to 
characteristic $\ell$--gerbes on $Z$ according to theorem \ref{6.9}. 

First of all, the above data determine a family pa\-ra\-me\-trized 
by $z \in Z$ of secondary characteristic classes of Godbillon--Vey type 
on the fibres $X_z$ of $\pi:Y \to Z$, according to theorem \ref{5.15}, 
namely 
\begin{equation}
\label{7.33}
[~\Lambda (C_1 Q, A_z)~] = [~\alpha_z \wedge Q (F_{A_z}^n)~], 
\end{equation}
where $\alpha = \beta \mid T (\pi)~, ~
d\alpha_z = C_1 (F_{A_z})$ and $F_A^{q+1}=0$. 
However, for $n > q$ and in particular for $\ell>0$, that is 
$\dim X=2n+1-\ell < 2n+1$, these forms on the fibres vanish 
identically. 

Next, we consider the case $n=q$, that is $\ell=0$ and 
$\dim X =2n+1=2q+1$ according to our general convention. 
Then the classes \eqref{7.33} actually live on the fibres 
$X_z = \pi^{-1} (z)$ and we obtain from 
theorem \ref{6.9} (i) a global $0$--gerbe 
\begin{equation}
\label{7.34}
[~\Lambda_{Y / Z} (C_1 Q, B)~] = \int_{Y/Z} ~
[~\beta \wedge Q (F_B^q)~] \in H_\D^1 (Z) = \Omega^0 (Z),  
\end{equation}
given fibrewise by 
\begin{equation}
\label{7.35}
[~\Lambda_{Y / Z} (C_1 Q, B)~] ~(z) = 
\int_{X_z} ~\alpha_z \wedge Q (F_{A_z}^q). 
\end{equation}
Thus the family of invariants in \eqref{7.35} are the integrated 
fibrewise Godbillon--Vey invariants, which are well--known to be 
variable and hence non--constant in $\Omega^0 (Z)$ for a suitable choice 
of the family of foliations (compare Heitsch~\cite{He1},~\cite{He2} and 
also the original work of Thurston~\cite{T}). 
A similar result is obtained for $n = q+\ell~, ~\ell>0$ and 
$\dim X=2n+1-\ell=2q+\ell+1$, 
in which case theorem \ref{6.9} (i) gives rise to (variable) 
characteristic $\ell$--gerbes 
\begin{equation}
\label{7.36}
[~\Lambda_{Y / Z} (C_1 Q, B)~] = \int_{Y/Z} ~
[~\beta \wedge Q (F_B^n)~] \in H_\D^{\ell+1} (Z) = 
\Omega^\ell (Z)/ d\Omega^{\ell-1} (Z),  
\end{equation}
determined by formula \eqref{7.31}; compare also \eqref{5.1b}. 

A more original class of gerbes is obtained in the 'rigid` range 
$n-\ell > q$, that is $\ell = 0, \ldots n-(q+1)$, in which case 
we still have $2q+1 < \dim X=2n+1-\ell$. 
Then we can invoke theorem \ref{6.9} (ii) to obtain well-defined 
flat characteristic Godbillon--Vey $\ell$--gerbes 
\begin{equation}
\label{7.37}
[~\Lambda_{Y / Z} (C_1 Q, B)~] = \int_{Y/Z} ~
[~\beta \wedge Q (F_B^n)~] \in H^{\ell} (Z, \R). 
\end{equation}
Note that for $n-\ell \geq q~, ~\ell > 0$, we have 
$2n+1 > \dim X=2n+1-\ell > 2q+1$. 
Hence, as already noted, the fibrewise classes vanish identically 
on the form level, while the forms $\Lambda (C_1 Q, B)$ are not 
necessarily closed on $Y$ unless $n \geq q + \dim Z$. 

In contrast, the classes investigated in Kotschick \cite{Kot}, 
Hoster--Kamber--Kotschick \cite{HKK}, are families of classes 
on a fixed manifold $X$, defined with respect to a $1$--parameter 
family $\F_t$ of foliations and foliated bundles and their suspension 
on the cylinder $X \x I$. Hoster in his thesis ~\cite{Ho} considers  
fibre spaces with flags of foliations along the fibres, but stays 
essentially in the context of ~\cite{HKK}. 

\end{example}



\end{document}